\documentclass[english,10pt]{article}

\usepackage[english]{babel}
\usepackage[utf8x]{inputenc}
\usepackage[T1]{fontenc}
\usepackage[formats]{listings}
\usepackage{geometry}
\usepackage{enumitem}

\usepackage{graphicx}
\usepackage[colorinlistoftodos]{todonotes}
\usepackage{setspace}
\usepackage{placeins}
\usepackage{enumitem}
\usepackage{titlesec}
\usepackage{cite}
\usepackage{comment}
\usepackage{caption}
\usepackage{subcaption}

\usepackage[colorlinks=true, allcolors=blue]{hyperref}

\usepackage{amsmath}
\usepackage{amsthm}
\usepackage{amssymb}
\usepackage{amsfonts}
\usepackage{bbm}
\usepackage{bm}
\usepackage{nicefrac}
\usepackage{mathtools}

\usepackage{tikz}

\usepackage{graphicx}
\usepackage{fullpage}
\usepackage{float}
\usepackage{wrapfig}
\usepackage{caption}

\usepackage{algpseudocode}
\usepackage{algorithm}
\usepackage{listings}

\usepackage{amsthm}\usepackage{dsfont}\usepackage{array}\usepackage{mathrsfs}\usepackage{cite}\usepackage{comment}\usepackage{mathrsfs}

\makeatother

\usepackage{babel}
\usepackage{color}
\usepackage{centernot}

\usepackage{babel}
\usepackage{sansmath}

\definecolor{revcolor}{RGB}{255,0,0}
\definecolor{ejc}{RGB}{139,0,139}
\definecolor{ps}{RGB}{0,0,200}

\title{A Modern Maximum-Likelihood Theory\\
  for High-dimensional Logistic Regression} 
\date{April 2018}

\author{Pragya Sur\thanks{Department of Statistics, Stanford
    University, Stanford, CA 94305, U.S.A.}  \and Emmanuel
  J. Cand\`es\footnotemark[1] \thanks{Department of Mathematics,
    Stanford University, Stanford, CA 94305, U.S.A.} }

\theoremstyle{plain}\newtheorem{theorem}{\textbf{Theorem}}

\theoremstyle{definition}

\theoremstyle{definition}

\newcommand{\bbeta}{\bm{\beta}}
\newcommand{\hbbeta}{\hat{\bm{\beta}}}
\newcommand{\hbeta}{\hat{\beta}}
\newcommand{\tbbeta}{\tilde{\bm{\beta}}}

\newcommand{\bX}{\bm{X}}

\newcommand{\bzero}{\bm{0}}

\newcommand{\dnorm}{{\mathcal{N}}}

\newcommand{\iid}{\stackrel{ \text{i.i.d.} }{\sim} }

\newcommand{\bu}{\bm{u}}

\newcommand{\R}{\mathbb{R}}

\newcommand{\eqd}{\stackrel{\mathrm{d}}{=}}

\newcommand{\taus}{\tau_{\star}}

\newcommand{\prox}{\mathsf{prox}}
\newcommand{\by}{\bm{y}}

\newcommand{\bS}{\bm{S}}

\newcommand{\bb}{\bm{b}}

\newcommand{\bZ}{\bm{Z}}
\newcommand{\bD}{\bm{D}}

\newcommand{\bH}{\bm{H}}

\newcommand{\tr}{\mathrm{Tr}}

\newcommand{\convd}{\stackrel{\mathrm{d}}{\rightarrow}}

\newcommand{\bU}{\bm{U}}
\newcommand{\bSigma}{\bm{\Sigma}}

\newcommand{\alphas}{\alpha_{\star}}

\newcommand{\prob}{\mathbb{P}}

\newcommand{\convP}{\stackrel{\prob}{\rightarrow}}

\newcommand{\bt}{\bm{t}}

\newcommand{\lambdas}{\lambda_{\star}}

\newcommand{\sigmas}{\sigma_{\star}}



\DeclareMathOperator{\E}{\mathbb{E}}

\begin{document}
\maketitle

\begin{abstract}
  Every student in statistics or data science learns early on that
  when the sample size $n$ largely exceeds the number $p$ of
  variables, fitting a logistic model produces estimates that are
  approximately unbiased. Every student also learns that there are
  formulas to predict the variability of these estimates which are
  used for the purpose of statistical inference; for instance, to
  produce p-values for testing the significance of regression
  coefficients. Although these formulas come from large sample
  asymptotics, 
  we are often told that we are on reasonably safe grounds when $n$ is
  large in such a way that $n \ge 5p$ or $n \ge 10p$. This paper shows
  that this is far from the case, and consequently, inferences
  routinely produced by common software packages are often unreliable.

  Consider a logistic model with independent features in which $n$ and
  $p$ become increasingly large in a fixed ratio. Then we show that
  (1) the MLE is biased, (2) the variability of the MLE is far greater
  than classically predicted, and (3) the commonly used
  likelihood-ratio test (LRT) is not distributed as a chi-square. The
  bias of the MLE is extremely problematic as it yields completely
  wrong predictions for the probability of a case based on observed
  values of the covariates.  We develop a new theory, which
  asymptotically predicts (1) the bias of the MLE, (2) the variability
  of the MLE, and (3) the distribution of the LRT. We empirically also
  demonstrate that these predictions are extremely accurate in finite
  samples. Further, an appealing feature is that these novel
  predictions depend on the unknown sequence of regression
  coefficients only through a single scalar, the overall strength of
  the signal. This suggests very concrete procedures to adjust
  inference; we describe one such procedure learning a single
  parameter from data and producing accurate inference.  For space
  reasons, we do not provide a full mathematical analysis of our
  results. However, we give a brief overview of the key arguments,
  which rely on the theory of (generalized) approximate message
  passing algorithms as well as on leave-one-observation/predictor out
  approaches.
\end{abstract}

\section{Introduction}
\subsection{Logistic regression: classical theory and practice} 

Logistic regression
\cite{lehmann2006testing,mccullagh1989generalized,van2000asymptotic}
is by and large the most frequently used model to estimate the
probability of a binary response from the value of multiple
features/predictor variables. It is used all the time in the social
sciences, the finance industry, the medical sciences, and so on. As an
example, a typical application of logistic regression may be to
predict the risk of developing a given coronary heart disease from a
patient's observed characteristics. Consequently, every graduate
student in statistics or any field that remotely involves data
analysis learns about logistic regression, perhaps before any other
nonlinear multivariate model. In particular, every student knows how
to interpret the excerpt of the computer output from Figure
\ref{fig:R}, which displays regression coefficient estimates, standard
errors and p-values for testing the significance of the regression
coefficients. In textbooks we learn the following:
\begin{enumerate}
\item Fitting a model via maximum likelihood produces estimates that
  are {\em approximately unbiased.}
\item There are formulas to {\em predict the accuracy or variability}
  of the maximum-likelihood estimate (MLE) (used in the computer
  output from Figure \ref{fig:R}).
\end{enumerate}
These approximations come from asymptotic results. Imagine we have $n$
independent observations $(y_i,\bX_i)$ where $y_i \in \{0,1 \}$ is the
response variable and $\bX_i \in \R^p$ the vector of predictor
variables. The logistic model posits that the probability of a case
conditional on the covariates is given by
\[
\mathbb{P}(y_i = 1 \, | \, \bX_i) = \rho'(\bX_i' \bbeta),  
\]
where $\rho'(t) = e^t/(1+e^{t})$ is the standard sigmoidal
function. Then everyone knows
\cite{lehmann2006testing,mccullagh1989generalized,van2000asymptotic}
that in the limit of large samples in which $p$ is fixed and
$n \rightarrow \infty$, the MLE $\hbbeta$ obeys
\begin{equation}\label{eq:asympnormclassical}
\sqrt{n}\bigl(\hbbeta - \bbeta \bigr) \convd
\dnorm(0,\bm{I}^{-1}_{\bbeta}),
\end{equation}
where $\bm{I}_{\bbeta}$ is the $p \times p$ Fisher information matrix
evaluated at the true $\bbeta$.  {A classical way of understanding
\eqref{eq:asympnormclassical} is in the case where the pairs
$(\bX_i, y_i)$ are i.i.d.~and the covariates $\bX_i$ are drawn
from a distribution obeying mild conditions so that the MLE exists and
is unique.}
Now the approximation \eqref{eq:asympnormclassical}
justifies the first claim of near unbiasedness. Further, software
packages then return standard errors by evaluating the inverse Fisher
information matrix at the MLE $\hbbeta$ (this is essentially what R
does in Figure \ref{fig:R}).  In turn, these standard errors are then
used for the purpose of statistical inference; for instance, they are
used to produce p-values for testing the significance of regression
coefficients, which researchers use in thousands of scientific
studies.
\begin{figure}
  \centering
\begin{minipage}{0.6\textwidth}
{\footnotesize
\begin{verbatim}
> fit = glm(y ~ X, family = binomial)
> summary(fit)

Coefficients:
            Estimate Std. Error z value Pr(>|z|)   
(Intercept)  0.25602    0.43191   0.593  0.55334   
X1           7.78102    4.09069   1.902  0.05715 . 
X2           9.80854    5.66019   1.733  0.08311 . 
X3          -8.14106    5.50490  -1.479  0.13917   
X4           0.01953    5.99945   0.003  0.99740   
X5          -5.18298    3.88752  -1.333  0.18245   
X6           9.48063    4.65335   2.037  0.04161 * 

...
---
Signif. codes:  0 ‘***’ 0.001 ‘**’ 0.01 ‘*’ 0.05 ‘.’ 0.1 ‘ ’ 1
\end{verbatim}
}
\end{minipage}
\caption{\small Excerpt from an object of class ``glm'' obtained by
  fitting a logistic model in R. The coefficient estimates
  $\hat \beta_j$ are obtained by maximum likelihood, and for each
  variable, R provides an estimate of the standard deviation of
  $\hat{\beta_j}$ as well as a p-value for testing whether
  $\beta_j = 0$ or not.}
\label{fig:R}
\end{figure}

Another well-known result in logistic regression is Wilks'
theorem \cite{wilks1938large}, which gives the asymptotic distribution
of the likelihood ratio test (LRT):
\begin{enumerate}
\setcounter{enumi}{2}
\item Consider the likelihood ratio obtained by dropping $k$ variables
  from the model under study. Then under the null hypothesis that none
  of the dropped variables belongs to the model, {\em twice the
    log-likelihood ratio (LLR) converges to a chi-square distribution} with
  $k$ degrees of freedom in the limit of large samples.
\end{enumerate}
Once more, this approximation is often used in lots of statistical
software packages to obtain p-values for testing the significance of
individual and/or groups of coefficients.

\subsection{Failures in moderately large dimensions} 

In modern-day data analysis, new technologies now produce extremely
large datasets, often with huge numbers of features on each of a
comparatively small number of experimental units. Yet, software
packages and practitioners continue to perform calculations as if
classical theory applies and, therefore, the main issue is this: do
these approximations hold in the modern setting where $p$ is not
vanishingly small compared to $n$?

To address this question, we begin by showing results from an
empirical study. Throughout this section, we set $n = 4000$ and unless
otherwise specified, $p = 800$ (so that the `dimensionality' $p/n$ is
equal to 1/5). We work with a matrix of covariates, which has
i.i.d.~$\dnorm(0,1/n)$ entries, and different types of regression
coefficients scaled in such a way that 
\[
\gamma^2 := \operatorname{Var}(\bX_i' \bbeta) =  5. 
\]
This is a crucial point: we want to make sure that the size of the
log-odds ratio $\bX_i' \bbeta$ does not increase with $n$ or $p$, so
that $\rho'(\bX_i'\bbeta)$ is not trivially equal to either $0$ or
$1$. Instead, we want to be in a regime where accurate estimates of
$\bbeta$ translate into a precise evaluation of a nontrivial
probability. With our scaling $\gamma = \sqrt{5} \approx 2.236$, about
95\% of the observations will be such that
$-4.472 \le \bX_i'\bbeta \le 4.472$ so that
$0.011 \le \rho'(\bX_i'\bbeta) \le 0.989$.



\begin{figure}
  \centering
  \includegraphics[scale=0.35,keepaspectratio]{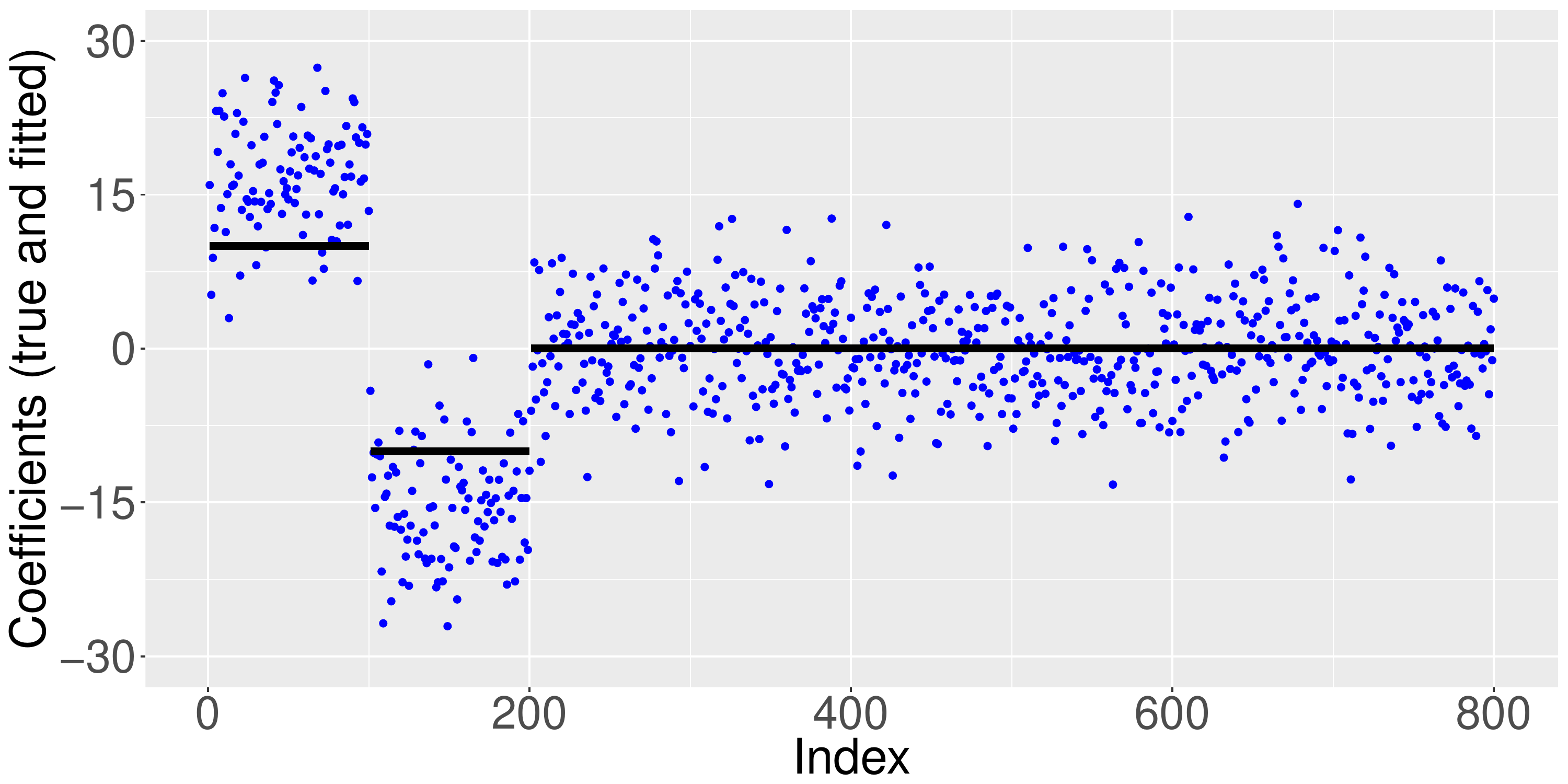}  
  \caption{True signal values $\beta_j$ in black and corresponding ML
    estimates $\hat{\beta}_j$ (blue points).  Observe that estimates
    of effect magnitudes are seriously biased upward.  }
\label{fig:centering_one}
\end{figure}

  
\paragraph{Unbiasedness?} 
Figure \ref{fig:centering_one} plots the true and fitted coefficients
in the setting where one quarter of the regression coefficients have a
magnitude equal to $10$, and the rest are zero. Half of the nonzero
coefficients are positive and the other half are negative.  A striking
feature is that the black curve does not pass through the center of
the blue scatter. This is in stark contradiction to what we would
expect from classical theory. Clearly, the regression estimates are
not close to being unbiased. When the true effect size $\beta_j$ is
positive, we see that the MLE has a strong tendency to overestimate
it.  Symmetrically, when $\beta_j$ is negative, the MLE tends to
underestimate the effect sizes in the sense that the fitted values are
in the same direction but with magnitudes that are too large. In other
words, for most indices $|\hat{\beta}_j| > |\beta_j|$ so that we are
over-estimating the magnitudes of the effects.

\begin{figure}
\begin{center}
	\begin{tabular}{ccc}
		\includegraphics[scale=0.3,keepaspectratio]{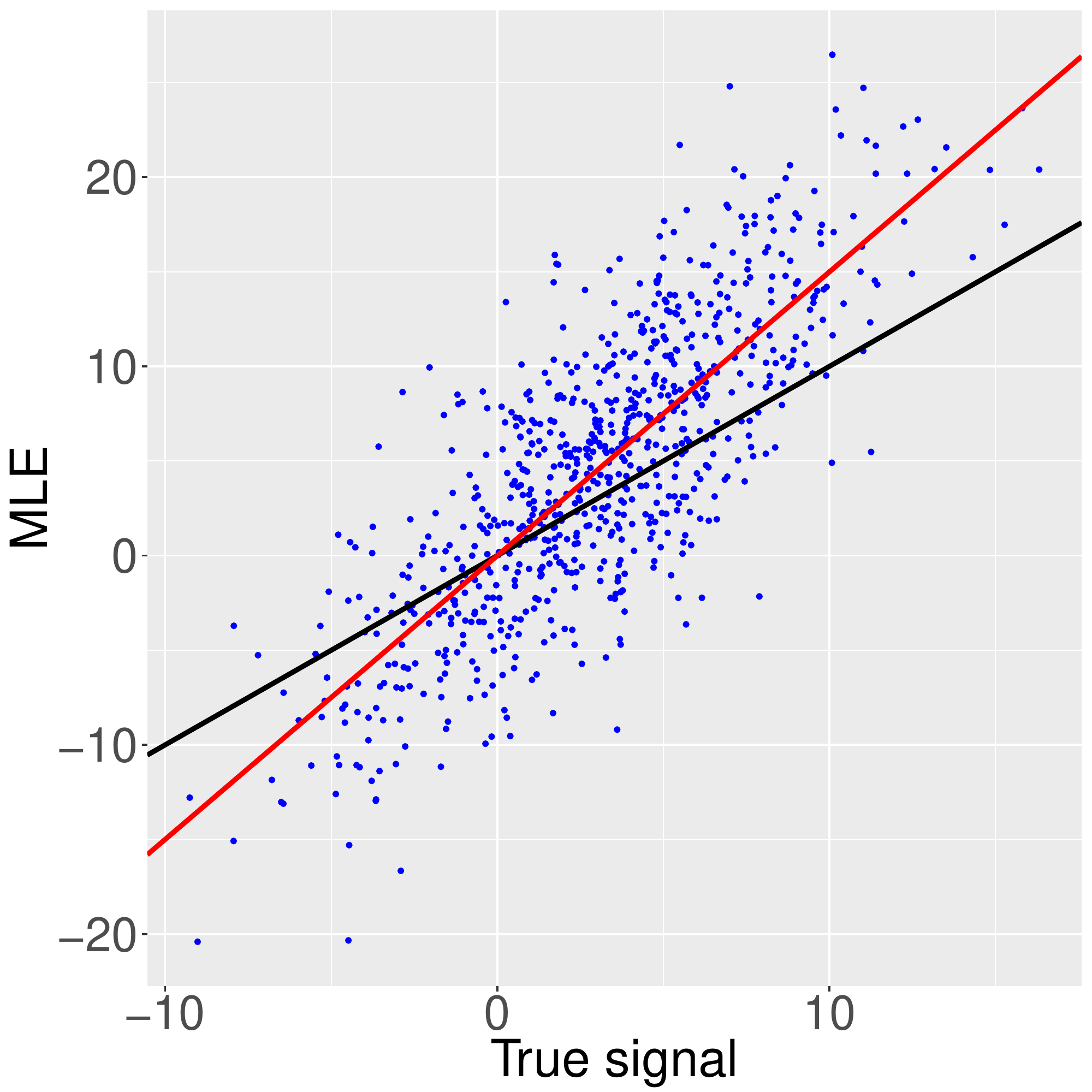} 
		& & \raisebox{0.4cm}{\includegraphics[scale=0.281,keepaspectratio]{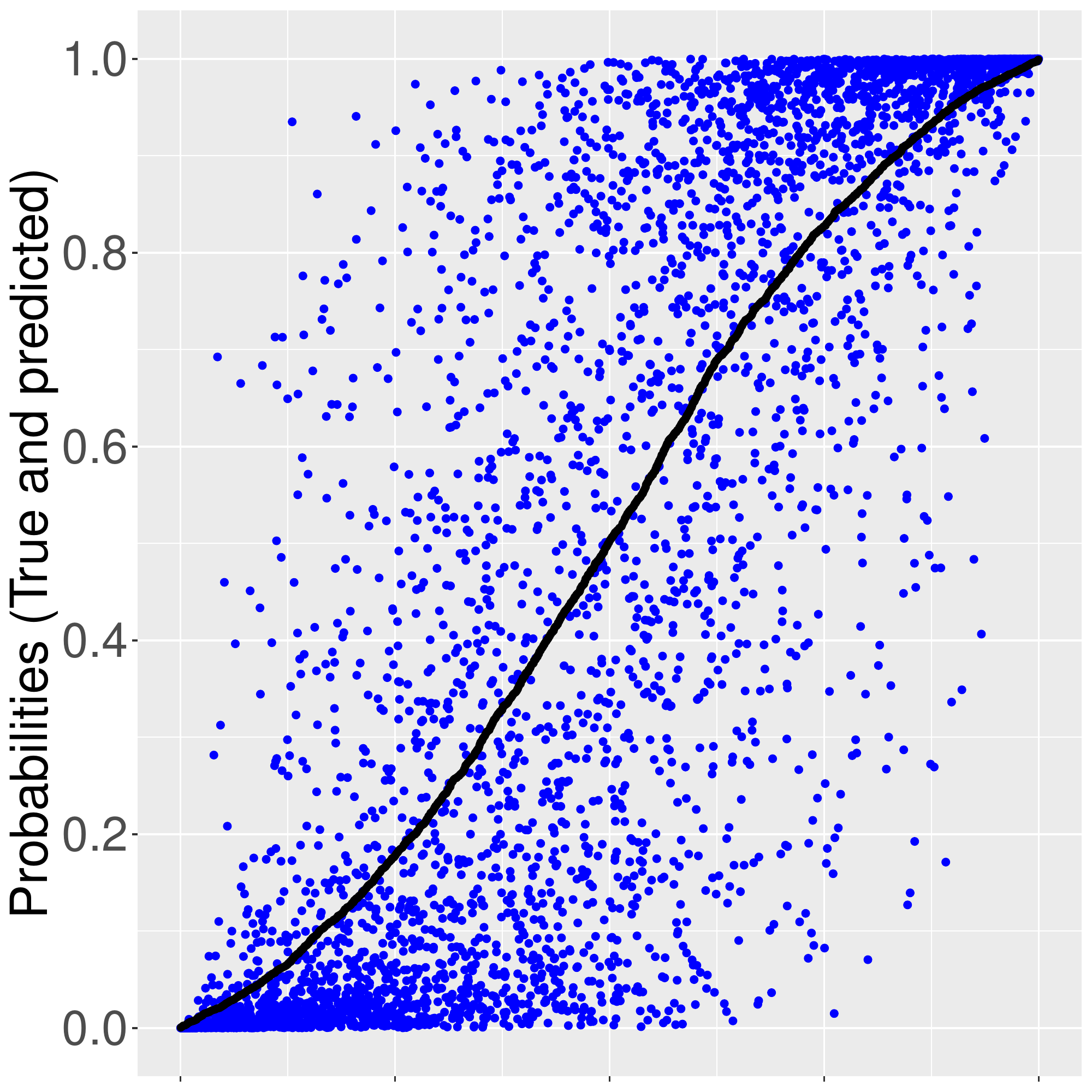}} \tabularnewline
		(a) & & (b)   \tabularnewline
	\end{tabular}
\end{center}
\caption{(a) Scatterplot of the pairs $(\beta_j, \hat{\beta}_j)$ for
  i.i.d. $\dnorm(3,16)$ regression coefficients. The black line has
  slope $1$. Again, we see that the MLE seriously overestimates effect
  magnitudes. The red line has slope $\alpha^\star \approx 1.499$
  predicted by the solution to \eqref{eq:main}. We can see that
  $\hat{\beta}_j$ seems centered around $\alpha^\star \beta_j$.  (b)
  True conditional probability $f(\bX_*) = \rho'(\bX_*' \bbeta)$
  (black curve), and corresponding estimated probabilities
  $\hat{f}(\bX_*) = \rho'(\bX_*' \hbbeta)$.  Observe the dramatic
  shrinkage of the estimates toward the end
  points.}\label{fig:centering}
\end{figure}

The strong bias is not specific to this example as the theory we
develop in this paper will make clear. Consider a case where the
entries of $\bbeta$ are drawn i.i.d.~from $\dnorm(3,16)$ (the setup is
otherwise unchanged). Figure \ref{fig:centering}(a), shows that the
pairs $(\beta_j, \hat{\beta}_j)$ are not distributed around a straight
line of slope 1; rather, they are distributed around a line with
a larger slope. Our theory predicts that the points should be scattered
around a line with slope $1.499$ shown in red, as if we could think
that $\E \hat{\beta}_j \approx 1.499 \beta_j$.

The strong bias is highly problematic for estimating the probability
of our binary response. Suppose we are given a vector of covariates
$\bX_*$ and estimate the regression function
$f(\bX_*) = \mathbb{P}(y = 1 \, | \, \bX_*)$ with
\[
\hat{f}(\bX_*) = \rho'(\bX_*'\hbbeta).
\]
Then because we tend to over-estimate the magnitudes of the effects,
we will also tend to over-estimate or under-estimate the probabilities
depending on whether $f(\bX_*)$ is greater or less than a half. This
is illustrated in Figure \ref{fig:centering}(b). Observe that when
$f(\bX_*) < 1/2$, lots of predictions tend to be close to zero, even
when $f(\bX_*)$ is nowhere near zero. A similar behavior is obtained
by symmetry; when $f(\bX_*) > 1/2$, we see a shrinkage toward the
other end point, namely, one. 
Hence, we see a massive shrinkage towards
the extremes and the phenomenon is amplified as the true probability
$f(\bX_*)$ approaches zero or one. Expressed differently, the MLE may
predict that an outcome is almost certain (i.e.~$\hat{f}$ is close to
zero or one) when, in fact, the outcome is not at all certain. This
behavior is misleading.

\begin{figure}
\begin{center}
	\begin{tabular}{ccc}
	\includegraphics[scale=0.3,keepaspectratio]{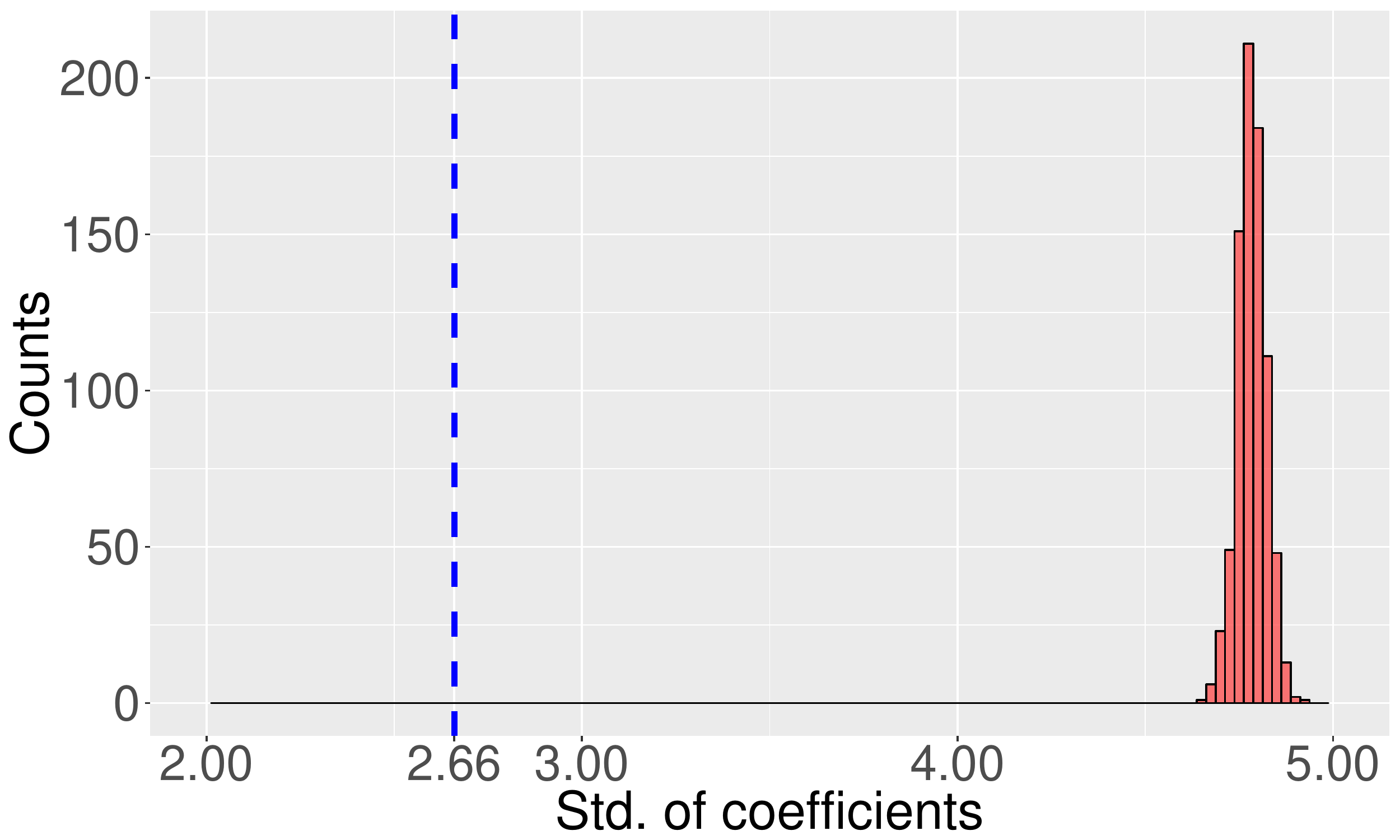} 
		& & \includegraphics[scale=0.3,keepaspectratio]{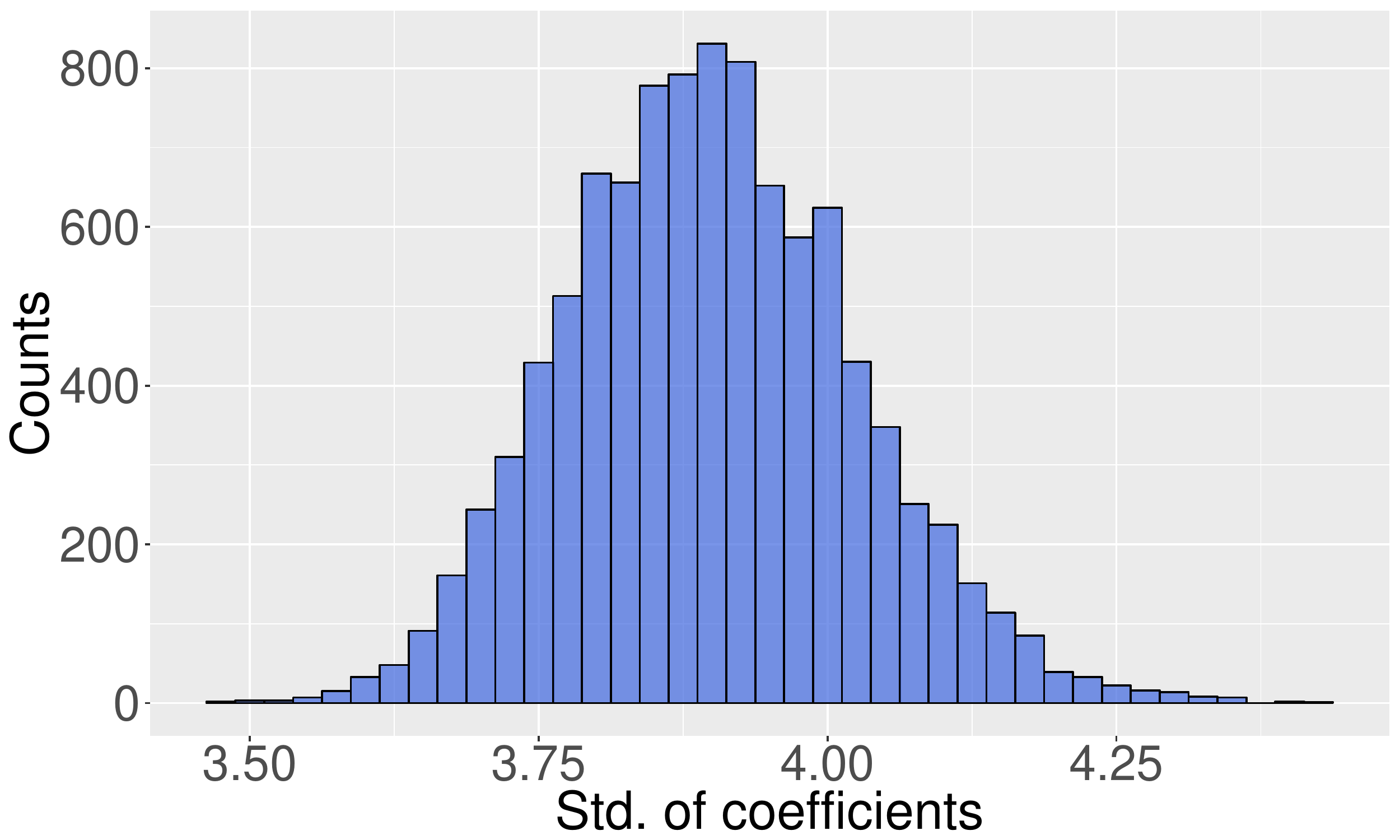}\label{fig:var} \tabularnewline
		(a) & & (b)  \tabularnewline
	\end{tabular}
\end{center}
\caption{(a) The histogram is the distribution of $\text{STD}(\hat{\beta}_j)$ for each variable $j$, in which the standard deviation is estimated from $10,000$ samples. The classical predicted standard error value is shown in blue.
  Classical theory underestimates the variability of the MLE. (b)
  Standard error estimates computed from R for a single null (for
  which $\beta_j = 0$) obtained across 10,000 replicates resampling
  the response vector and the covariate matrix.}
\label{fig:varandpval}
\end{figure}

\paragraph{Accuracy of classical standard errors?} 
Consider the same matrix $\bX$ as before and regression coefficients
now sampled as follows: half of the $\beta_j$'s are i.i.d.~draws from
$\dnorm(7,1)$, and the other half vanish. Figure
\ref{fig:varandpval}(a) shows standard errors computed via Monte Carlo
of ML estimates $\hat{\beta}_j$ corresponding to null
coordinates. This is obtained by fixing the signal $\bbeta$ and
resampling the response vector and covariate matrix $10,000$ times.
Note that for any null coordinate, the classical prediction for the
standard deviation based on the inverse Fisher information can be
explicitly calculated in this setting and turns out to be equal to
$2.66$, see Appendix \ref{sec:Fisher}. Since the standard deviation
values evidently concentrate around $4.75$, we see that in higher
dimensions, the variance of the MLE is likely to be much larger than
that predicted classically.  Naturally, using classical predictions
would lead to grossly incorrect p-values and confidence statements, a
major issue first noticed in \cite{candes2016panning}.

The variance estimates obtained from
statistical software packages
are different
from the value 2.66 above because they do not take expectation over
the covariates and use the MLE $\hbbeta$ in lieu of $\bbeta$ (plugin
estimate), see Appendix \ref{sec:Fisher}. 
Since practitioners use these estimates all the time, it is useful to
describe how they behave.
To this end, for each of the $10,000$
samples $(\bX, \by)$ drawn above, we obtain the R standard error
estimate for a single MLE coordinate corresponding to a null
variable. The histogram is shown in Figure \ref{fig:varandpval}
(b). The behavior for this specific coordinate is typical of that
observed for any other null coordinate, and the maximum value for
these standard errors remains below $4.5$, significantly below the
typical values observed via Monte Carlo simulations in Figure
\ref{fig:varandpval}(a).  


\paragraph{Distribution of the LRT?} By now, the reader should be
suspicious that the chi-square approximation for the distribution of
the likelihood-ratio test holds in higher dimensions. Indeed, it does
not and this actually is not a new observation.  In
\cite{sur2017likelihood}, the authors established that for a class of
logistic regression models, the LRT converges weakly to a {\em
  multiple} of a chi-square variable in an asymptotic regime in which
both $n$ and $p$ tend to infinity in such a way that
$p/n \rightarrow \kappa \in (0,1/2)$.  The multiplicative factor is an
increasing function of the limiting aspect ratio $\kappa$, and exceeds
one as soon as $\kappa$ is positive. This factor can be computed by
solving a nonlinear system of two equations in two unknowns given in
\eqref{eq:reduced} below. Furthermore, \cite{sur2017likelihood} links
the distribution of the LRT with the asymptotic variance of the
marginals of the MLE, which turns out to be provably higher than that
given by the inverse Fisher information. These findings are of course
completely in line with the conclusions from the previous
paragraphs. The issue is that the results from
\cite{sur2017likelihood} assume that $\bbeta = 0$; that is, they apply
under the global null where the response does not depend upon the
predictor variables, and it is a priori not clear how the theory would
extend beyond this case. Our goal in this paper is to study properties
of the MLE and the LRT for high-dimensional logistic regression models
under general signal strengths---restricting to the regime where the
MLE exists.

\begin{figure}
  \centering
  \includegraphics[scale=0.35,keepaspectratio]{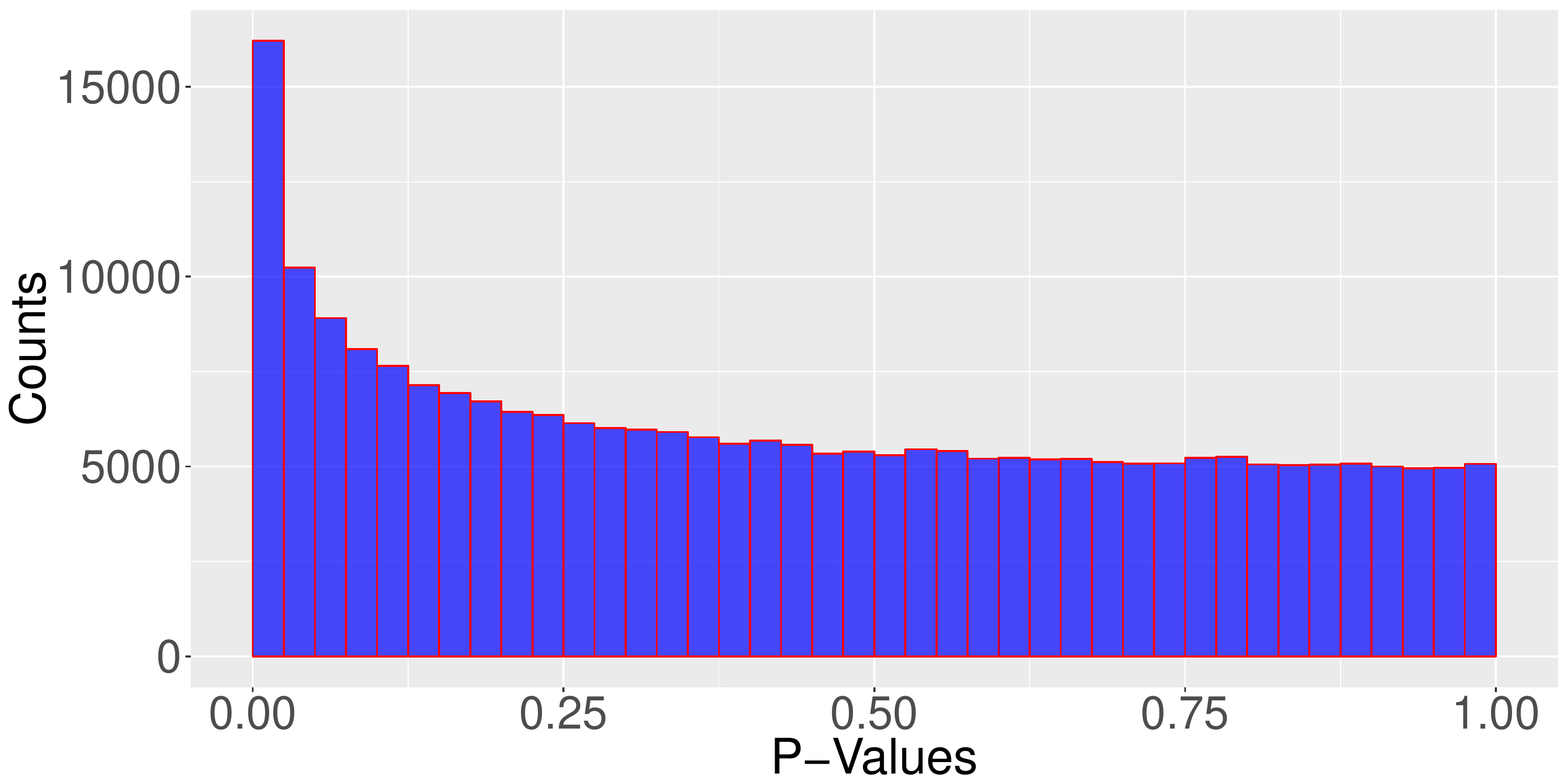}  
  \caption{ P-values calculated
  from the $\chi^2_1$ approximation to the LLR. Parameters:
  $n=4000, \kappa =0.2$, with half the co-ordinates of $\bbeta$
  non-zero, generated i.i.d. from
  $\dnorm(7,1)$. }
\label{fig:LRT}
\end{figure}  
To investigate what happens when we are not under the global null,
consider the same setting as in Figure \ref{fig:varandpval}. Figure
\ref{fig:LRT} shows the histogram of the p-values for testing a null
coefficient based on the chi-square approximation.  Not only are the
p-values far from uniform, the enormous mass near zero is extremely
problematic for multiple testing applications, where one examines
p-values at very high levels of significance, e.g.~near Bonferroni
levels. In such applications, one would be bound to make a very large
number of false discoveries from using p-values produced by software
packages. To further demonstrate the large inflation near the small
p-values, we display in Table \ref{tab:pvalfits} estimates of the
p-value probabilities in bins near zero. The estimates are much higher
than what is expected from a uniform distribution. Clearly, the
distribution of the LRT is far from a $\chi^2_1$.
\begin{table}[h!]
\centering
    \begin{tabular} {|c|c|}
    \hline
     & Classical \\ \hline
    $\prob\{\text{p-value}  \leq 5 \% \}$ & $10.77 \% (0.062 \%)$ \\ 
    \hline
        $\prob\{\text{p-value} \leq 1 \% \}$ & $3.34\% (0.036 \%)$   \\ 
        \hline
        $\prob\{\text{p-value}  \leq 0.5 \% \}$& $1.98 \% (0.028 \%)$  \\ 
        \hline
        $\prob\{\text{p-value}  \leq 0.1 \% \}$ &  $0.627 \% (0.016 \%)$ \\ \hline
        $\prob\{\text{p-value}  \leq 0.05 \% \}$ &  $0.365 \% (0.012 \%)$  \\ \hline
        $\prob\{\text{p-value}  \leq 0.01 \% \}$ &  $0.136 \% (0.007 \%)$  \\ \hline
    \end{tabular} \\ 
    \caption{P-value probabilities with standard errors in
      parentheses. Here, $n=4000$, $p = 800$, $\bX$ has
      i.i.d.~Gaussian entries, and half of the entries of $\bbeta$ are
      drawn from $\dnorm(7,1)$.}
    \label{tab:pvalfits}   
\end{table}

\paragraph{Summary.} We have hopefully made the case that classical
results, which software packages continue to rely upon, are downright
erroneous in higher dimensions.
\begin{enumerate}
\item Estimates seem systematically biased in the sense that effect magnitudes are overestimated.
\item Estimates are far more variable than classically predicted.
\item Inference measures, e.g. p-values, are unreliable especially at
  small values.
\end{enumerate}
Given the widespread use of logistic regression in high dimensions, a
novel theory explaining how to adjust inference to make it valid is
seriously needed.

\subsection{Our contribution}

Our contribution is to develop a brand new theory, which applies to
high-dimensional logistic regression models with independent
variables, and is capable of accurately describing all the phenomena
we have discussed. Taking them one by one, the theory from this paper
predicts:
\begin{enumerate}
\item the bias of the MLE;
\item the variability of the MLE;
\item and the distribution of the LRT.
\end{enumerate}
These predictions are, in fact, asymptotically exact in a regime where
the sample size and the number of features grow to infinity in a fixed
ratio.  Moreover, we shall see that our theoretical predictions are
extremely accurate in finite sample settings in which $p$ is a
fraction of $n$, e.g.~$p = 0.2 n$.

A very useful feature of this novel theory is that in our model, all
of our predictions depend on the true coefficients $\bbeta$ only
through the signal strength $\gamma$, where
$\gamma^2 := \operatorname{Var}(\bX_i'\bbeta)$.  This immediately
suggests that estimating some high-dimensional parameter is not
required to adjust inference. We propose in Section \ref{sec:gammaest}
a method for estimating $\gamma$ and empirically study the quality of
inference based on this estimate.

At the mathematical level, our arguments are very involved. Our
strategy is to introduce an approximate message passing algorithm that
tracks the MLE in the limit of a large number of features and
samples. In truth, a careful mathematical analysis is delicate and
requires a great number of steps. This is why in this expository paper
we have decided to provide the reader only with the main ideas. All
the details may be found in the separate document \cite{sur2018modernproofs}.

\subsection{Prior work} 

Asymptotic properties of M-estimators in the context of linear
regression have been extensively studied in diverging dimensions
starting from \cite{huber1973robust}, followed by
\cite{portnoy1984asymptotic} and \cite{portnoy1985asymptotic}. These
papers investigated the consistency and asymptotic normality
properties of M-estimators in a regime where $p = o(n^{\alpha})$, for
some $\alpha < 1$. Later on, the regime where $p$ is comparable to $n$
became the subject of a series of remarkable works \cite{el2013robust,
  bean2013optimal, donoho2013high, karoui2013asymptotic,
  el2015impact}; these works only concern the linear model. The
finding in these papers is that although M-estimators remain
asymptotically unbiased, they are shown to exhibit a form of `variance
inflation'.

Moving on to more general exponential families,
\cite{portnoy1988asymptotic} studied the asymptotic behavior of
likelihood methods and established that the classical Wilks' theorem
holds if $p^{3/2}/n \rightarrow 0$ and, moreover, that the classical
normal approximation to the MLE holds if $p^2/n \rightarrow
0$. Subsequently, \cite{he2000parameters} quantified the $\ell_2$
estimation error of the MLE when $p \log p /n \rightarrow 0$. Very
recently, the authors from \cite{fan2017nonuniform} investigated the
classical asymptotic normality of the MLE under the global null and
regimes in which it may break down as the dimensionality increases.
In parallel, there also exists an extensive body of literature on
penalized maximum likelihood estimates/procedures for generalized
linear models, see \cite{van2008high, kakade2010learning,
  bunea2008honest, jankova2017honest, van2014asymptotically,
  fan2011nonconcave,belloni2013honest}, for example. This body of
literature often allows $p$ to be larger than $n$ but relies upon
strong assumptions about the extreme sparsity of the underlying
signal. The setting in these works is, therefore, completely different
from ours.

Finite sample behavior of both the MLE and the LRT have been
extensively studied in the literature. It has been observed that when
the sample size is small, the MLE is found to be biased for the
regression coefficients. In this context, a series of works---
\cite{firth1993bias,anderson1979logistic,mclachlan1980note,schaefer1983bias,copas1988binary,cordeiro1991bias,cordeiro2014introduction}
and the references therein---proposed finite sample corrections to
the MLE, which typically hinges on an asymptotic expansion of the MLE
up to $O(1/n)$ terms. One can plug in an estimator of the $O(1/n)$
term, which would make the resultant corrected statistic $o(1/n)$
accurate. All of these works are in the low-dimensional setting where
the MLE is still asymptotically unbiased. The observed bias was simply
attributed to a finite sample effect. Jackknife bias reduction
procedures for finite samples have been proposed (see
\cite{bull1994two} and the references cited therein for other finite
sample corrections).  Similar corrections for the LRT have been
studied, see for instance,
\cite{bartlett1937properties,box1949general,lawley1956general,cordeiro1983improved,
  cordeiro1995bartlett,bickel1990decomposition,cribari1996bartlett,cordeiro2014introduction,moulton1993bartlett}. It
was demonstrated in \cite{sur2017likelihood} that such finite sample
corrections do not yield valid p-values in the high-dimensional regime
we consider. In \cite{jennings1986judging}, the author proposed a
measure for detecting inadequacy of inference that was based on
explicit computation of the third order term in the Taylor expansion
of the likelihood function. This term is known to be asymptotically
negligible in the low-dimensional setting, and is found to be
negligible asymptotically in our high-dimensional regime as well, as
will be shown in this paper. Thus, this proposal also falls under the
niche of a finite sample correction.

A line of simulation based results exist to guide practitioners
regarding how large sample sizes are needed so that such finite sample
problems would not arise while using classical inference for logistic
regression. The rule of thumb is usually $10$ events per variable
(EPV) or more as mentioned in
\cite{peduzzi1996simulation,hosmer2013applied}, while a later study
\cite{vittinghoff2007relaxing} suggested that it could be even
less. As we clearly see in this paper, such a rule is not at all valid
when the number of features is large.
\cite{courvoisier2011performance} contested the previously established
$10$ EPV rule.

  
To the best of our knowledge, logistic regression in the regime where
$p$ is comparable to $n$ has been quite sparsely studied. As already
mentioned, this paper follows up on the earlier contribution
\cite{sur2017likelihood} of the authors, which characterized the LLR
distribution in the case where there is no signal (global null). This
earlier reference derived the asymptotic distribution of the LLR as a
function of the limiting ratio $p/n$. As we will see later in Section
\ref{sec:main}, this former result may be seen as a special case of
the novel Theorem \ref{thm:lrt}, which deals with general signal
strengths. As is expected, the arguments are now more complicated than
when working under the global null.

\section{Main Results}
\label{sec:main}

\paragraph{Setting.} We describe the asymptotic properties of the MLE
and the LRT in a high-dimensional regime, where $n$ and $p$ both go to
infinity in such a way that $p/n \rightarrow \kappa$. We work with
independent observations $\{\bX_i, y_i\}$ from a logistic model such
that $\mathbb{P}(y_i = 1 \, | \, \bX_i) = \rho'(\bX_i'\bbeta)$. We
assume here that $\bX_i \sim \dnorm(\bzero, n^{-1} \bm{I}_{p})$, where
$\bm{I}_p$ is the $p$-dimensional identity matrix. (This means that
the columns of the matrix $\bX$ of covariates are unit-normed in the
limit of large samples.). The exact scaling of $\bX_i$ is not
important. As noted before, the important scaling is the signal
strength $\bX_i'\bbeta$ and we assume that the $p$ regression
coefficients (recall that $p$ increases with $n$) are scaled in such a
way that
\begin{equation}
  \label{eq:scaling}
\lim_{n \rightarrow \infty}  \, \operatorname{Var}(\bX_i' \bbeta) = \gamma^2,
\end{equation}
where $\gamma$ is fixed.  It is useful to think of the parameter
$\gamma$ as the signal strength.  Another way to express
\eqref{eq:scaling} is to say that
$\lim_{n \rightarrow \infty} \,\|\bbeta\|^2/n = \gamma^2$.

\subsection{When does the MLE exist?} 

\begin{figure}
\begin{center}
	\begin{tabular}{cc}
	\includegraphics[scale=0.3,keepaspectratio]{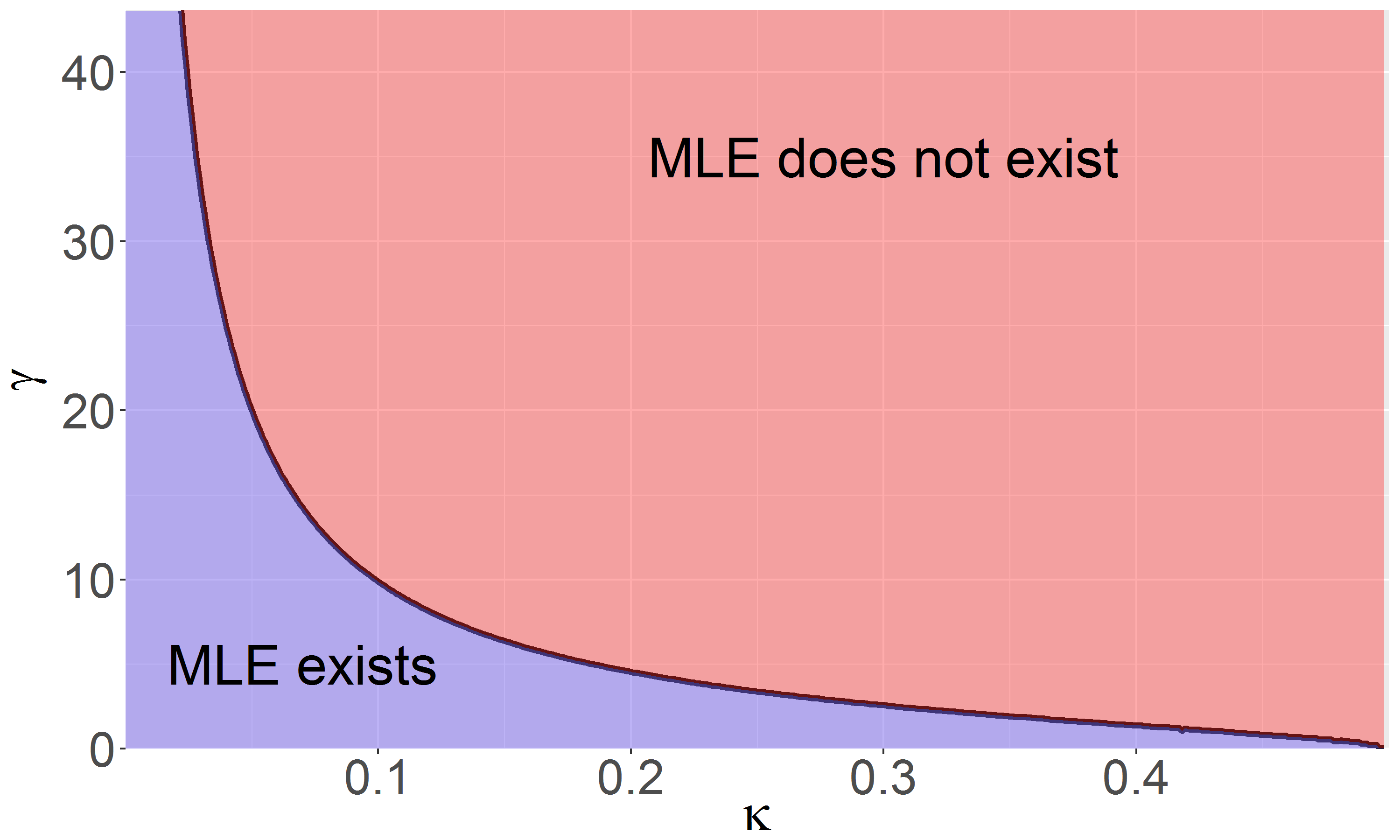} 
		&  \includegraphics[scale=0.3,keepaspectratio]{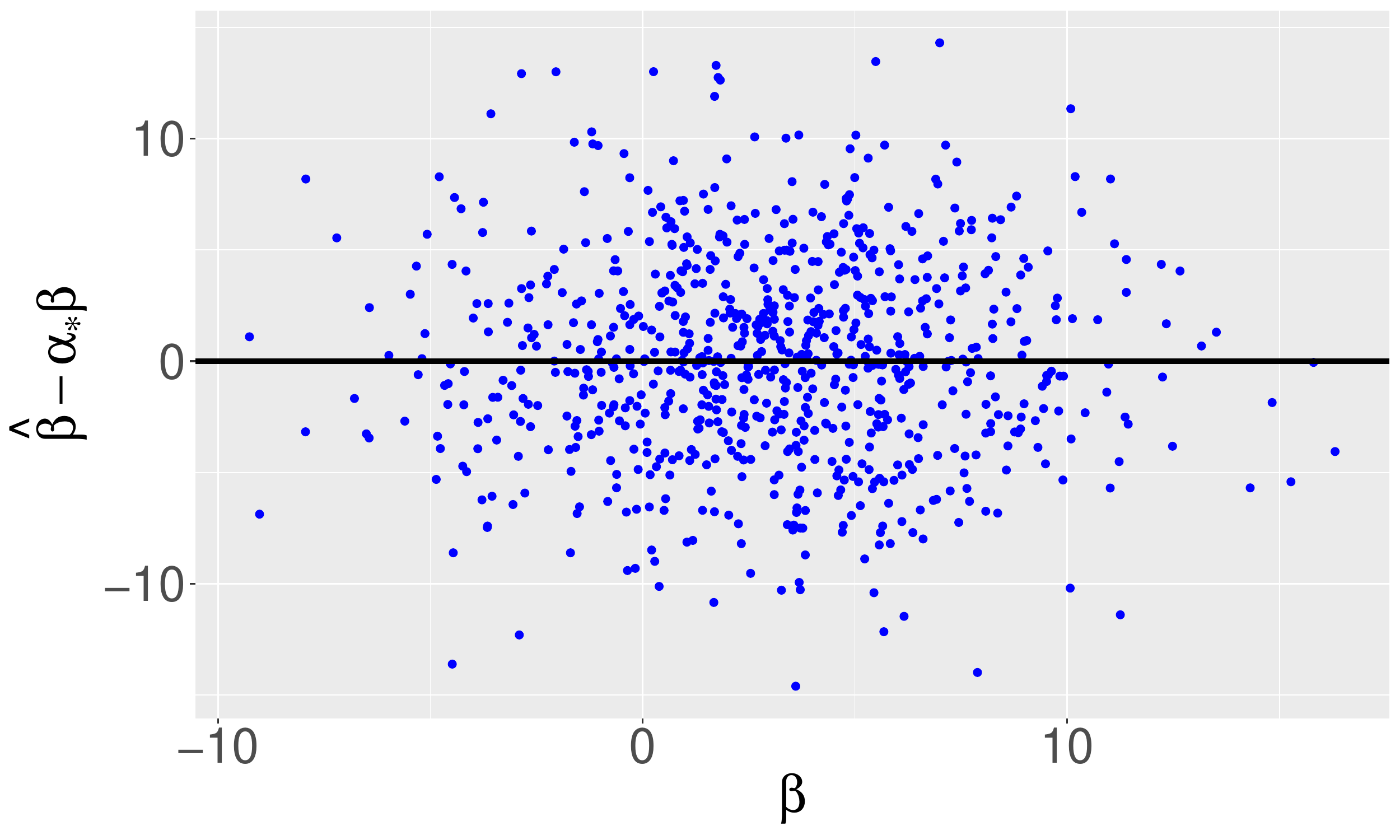}\label{fig:var} \tabularnewline
		(a) & (b)  \tabularnewline
	\end{tabular}
\end{center}
\caption{(a) Regions in which the MLE asymptotically exists and is
  unique and that in which it does not. The boundary curve is
  explicit and given by \eqref{eq:gmle}. (b) In the setting of Figure
  \ref{fig:centering}, scatterplot of the centered MLE
  $\hbeta_j - \alphas \beta_j$ vs.~the true signal $\beta_j$.}
\label{fig:kappagamma}
\end{figure}

The MLE $\hbbeta$ is the minimizer of the negative log-likelihood
$\ell$ defined via (observe that the sigmoid is the first derivative
of $\rho$)
\begin{equation}\label{eq:likelhood}
  \ell(\bm{b})=\sum_{i=1}^n \{ \rho(\bX_i'\bm{b}) - y_i \, (\bX_i'\bm{b}) \}, \qquad \rho(t) = \log (1+e^t).
 \end{equation}
 A first important remark is that in high dimensions, the MLE does not
 asymptotically exist if the signal strength $\gamma$ exceeds a
 certain functional $g_{\text{MLE}}(\kappa)$ of the dimensionality:
 i.e.~$\gamma > g_{\text{MLE}}(\kappa)$.  This happens because in such
 cases, there is a perfect separating hyperplane---separating the
 cases from the controls if you will---sending the MLE to infinity.
 In \cite{sur2017likelihood}, the authors proved that if $\gamma = 0$
 then $g_{\text{MLE}}(1/2) = 0$ (to be exact, they assumed
 $\bbeta = \bzero$). To be more precise, the MLE exists if
 $\kappa < 1/2$ whereas it does not if $\kappa > 1/2$
 \cite{cover1965geometrical, cover1964thesis}. Here, it turns out that
 a companion paper \cite{sur2018mle} precisely characterizes the
 region in which the MLE exists.
 \begin{theorem}[\cite{sur2018mle}]  
   \label{thm:mle_exist} Let $Z$ be a standard normal variable with
   density $\varphi(t)$ and $V$ be an independent continuous random
   variable with density $2\rho'(\gamma t) \varphi(t)$. With
   $x_+ = \max(x,0)$, set
   \begin{equation}
     \label{eq:gmle}
     g^{-1}_{\text{\em MLE}}(\gamma) = \min_{t \in \R} \, \left\{\E (Z -
       tV)_+^2\right\},  
   \end{equation}
   which is a decreasing function of $\gamma$. Then in the setting
   described above,
  \[
\begin{array}{lll}
  \gamma > g_{\text{\em MLE}}(\kappa) 
  & \quad \implies \quad 
  &
    \lim_{n,p
    \rightarrow
    \infty} 
    \mathbb{P}\{\text{\em MLE
    exists}\}
    \rightarrow
    0, \\
  \gamma < g_{\text{\em MLE}}(\kappa)  
  & \quad \implies \quad 
  &
    \lim_{n,p
    \rightarrow
    \infty} 
    \mathbb{P}\{\text{\em MLE
    exists}\}
    \rightarrow
    1.
\end{array}                     
\]
\end{theorem}
Hence, the curve $\gamma = g_{\text{MLE}}(\kappa)$, or, equivalently,
$\kappa = g^{-1}_{\text{MLE}}(\gamma)$ shown in Figure
\ref{fig:kappagamma} separates the $\kappa$--$\gamma$ plane into two
regions: one in which the MLE asymptotically exists and one in which
it does not. Clearly, we are interested in this paper in the former
region (the purple region in Figure \ref{fig:kappagamma}).


\subsection{A system of nonlinear equations} 
As we shall soon see, the asymptotic behavior of both the MLE and the
LRT is characterized by a system of equations in three variables
$(\alpha,\sigma,\lambda)$:
\begin{equation}
  \left\{ \begin{aligned} 
     \sigma^2  & = \frac{1}{\kappa^2} \E  \left[2\rho'(Q_1)\left(\lambda \rho'(\prox_{\lambda \rho}(Q_2))\right)^2 \right]  \\
     0 & = \E \left[\rho'(Q_1)Q_1 \lambda \rho'(\prox_{\lambda \rho}(Q_2))\right] \\
    1- \kappa  &  = \E \left[\frac{2\rho'(Q_1)}{1+\lambda \rho''(\prox_{\lambda \rho}( Q_2))} \right]
  \end{aligned}
\right.
  \label{eq:main}
\end{equation}
where $(Q_1,Q_2)$ is a bivariate normal variable with mean $\bzero$
and covariance 
\begin{equation}\label{eq:covfunc}
\bSigma(\alpha,\sigma) = \left[ {\begin{array}{cc}
  \gamma^2 &  -\alpha \gamma^2 \\
   -\alpha \gamma^2 & \alpha^2 \gamma^2 + \kappa \sigma^2 \end{array} } \right] .
\end{equation}
With $\rho$ as in \eqref{eq:likelhood}, the proximal mapping operator
is defined via
\begin{align}\label{eq:prox}
  \prox_{\lambda \rho}(z)  = \arg \min_{t \in \R} \left\{ \lambda \rho(t) + \frac{1}{2} (t -z )^2 \right\}.
\end{align}

The system of equations \eqref{eq:main} is parameterized by the pair
$(\kappa,\gamma)$ of dimensionality and signal strength parameters. It
turns out that the system admits a unique solution if and only if
$(\kappa,\gamma)$ is in the region where the MLE asymptotically
exists!


It is instructive to note that in the case where the signal strength
vanishes, $\gamma=0$, the system of equations \eqref{eq:main} reduces
to the following two-dimensional system:
\begin{equation}
  \left\{ \begin{aligned} 
     \sigma^2  & = \frac{1}{\kappa^2} \E  \left[\left(\lambda \rho'(\prox_{\lambda \rho}(\tau Z))\right)^2 \right]  \\
    1- \kappa  &  = \E \left[\frac{1}{1+\lambda \rho''(\prox_{\lambda \rho}( \tau Z))} \right] 
  \end{aligned}
\right.  \quad \tau^2 := \kappa \sigma^2, \quad Z \sim \dnorm(0,1). 
\label{eq:reduced}
\end{equation}
This holds because $Q_1 = 0$. It is not surprising that this system be
that from \cite{sur2017likelihood} since that work considers
$\bbeta = 0$ and, therefore, $\gamma = 0$.
\begin{figure}
\begin{center}
	\begin{tabular}{ccc}
\includegraphics[scale=0.3,keepaspectratio]{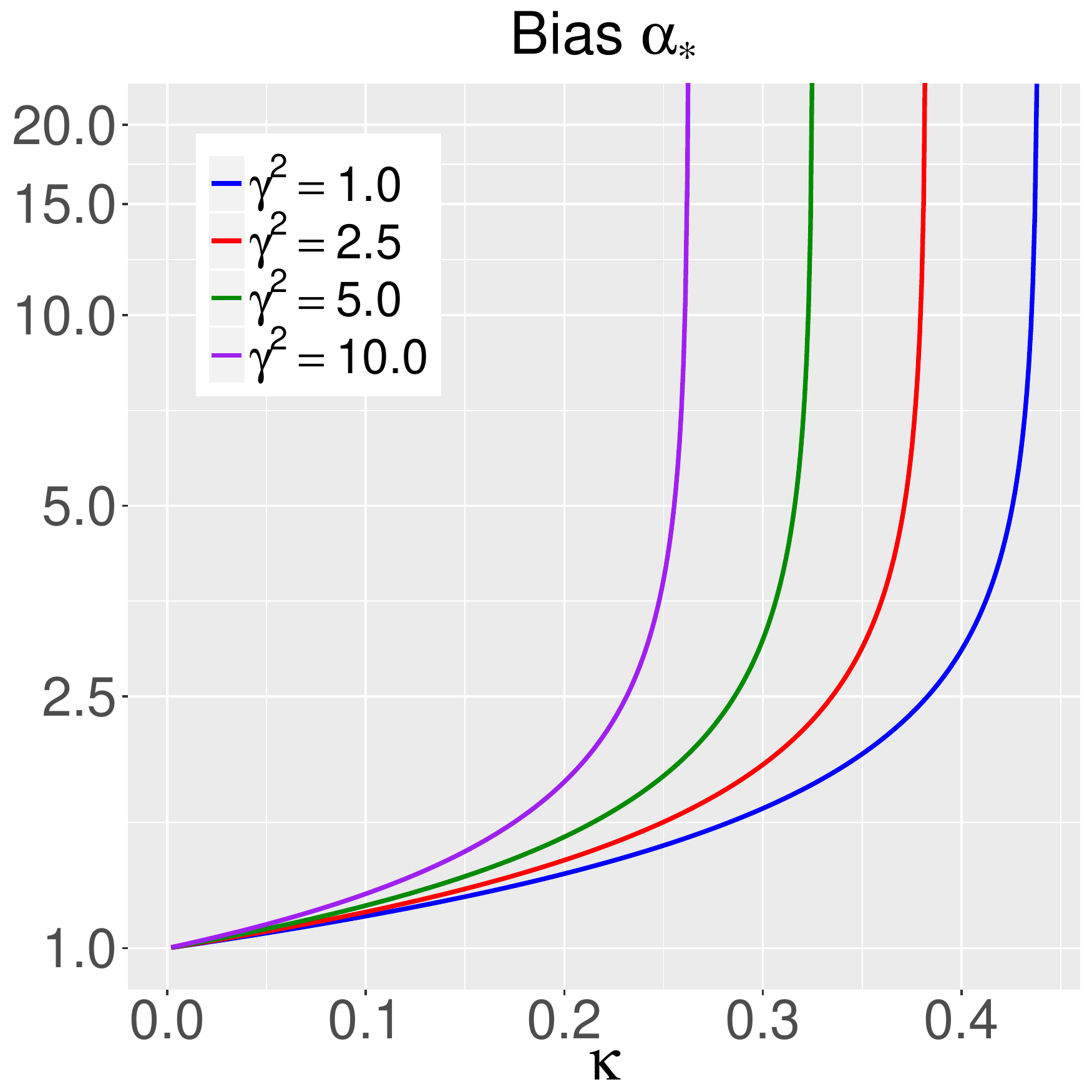} 
& \includegraphics[scale=0.3,keepaspectratio]{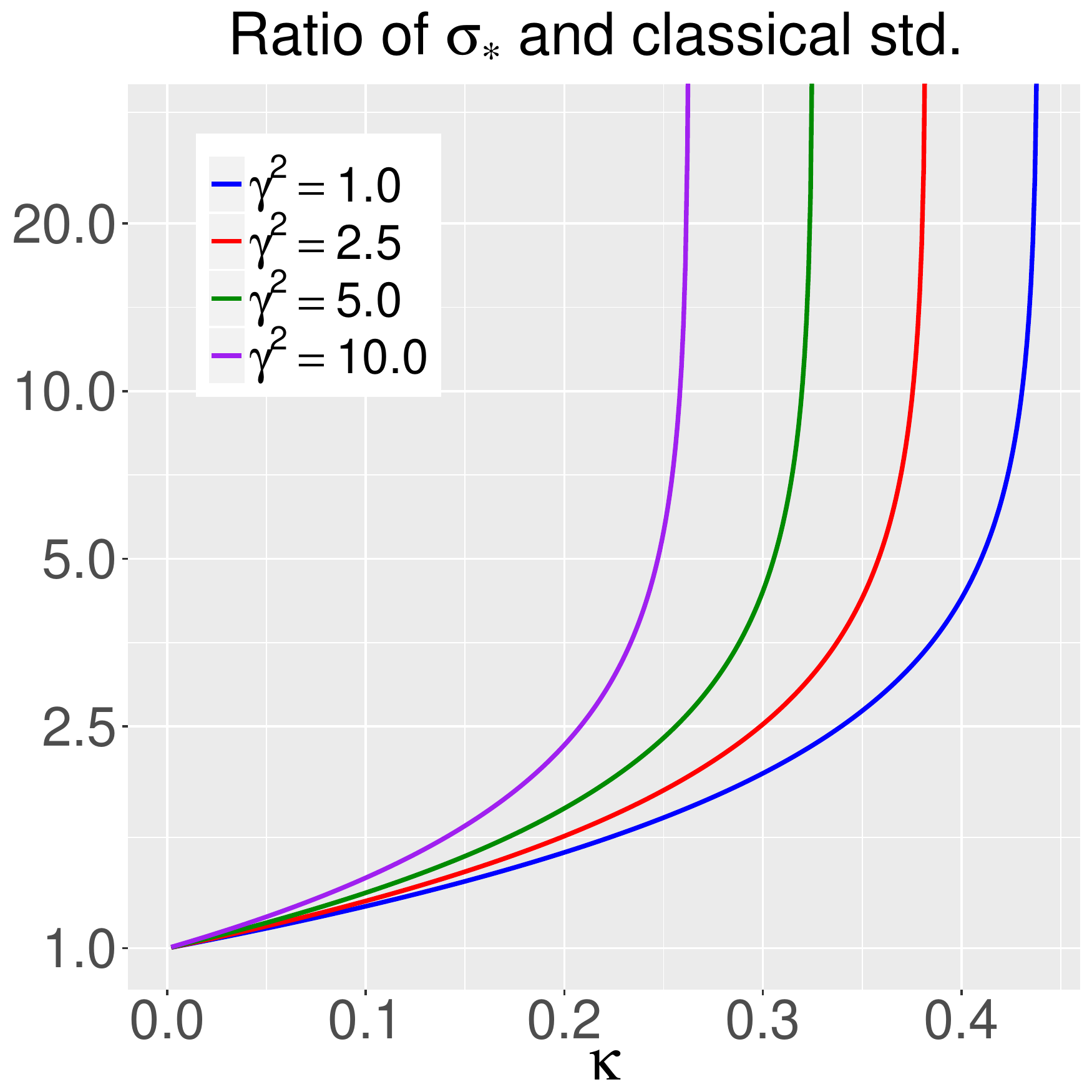} 
& \includegraphics[scale=0.3,keepaspectratio]{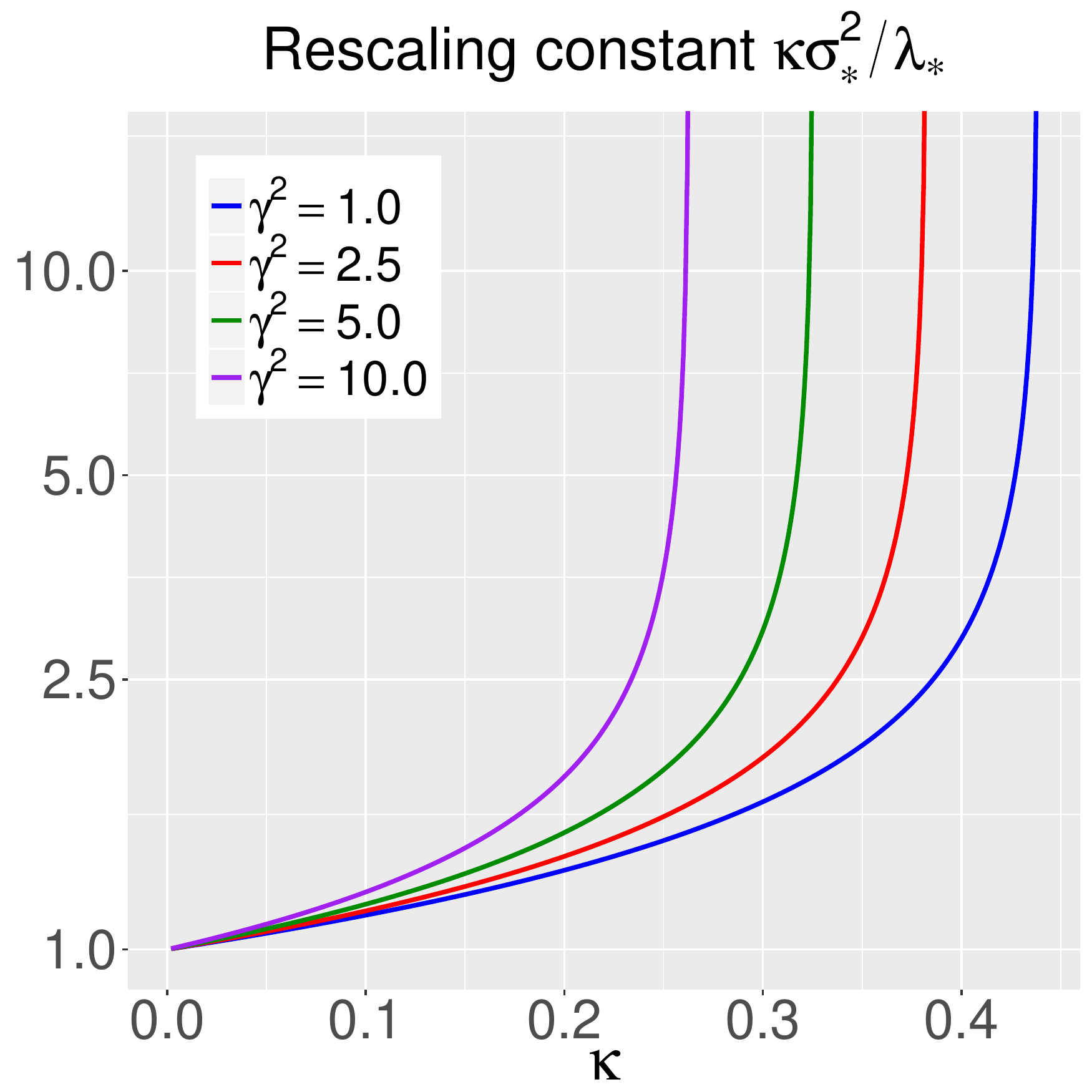}\tabularnewline
		(a) & (b) & (c)  \tabularnewline
	\end{tabular}
\end{center}
\caption{ (a) Bias $\alphas$ as a function of $\kappa$, for different
  values of the signal strength $\gamma$. Note the logarithmic scale
  for the y-axis. The curves asymptote at the value of $\kappa$ for
  which the MLE ceases to exist. (b) Ratio of the theoretical
  prediction $\sigmas$ and the average standard deviation of the
  coordinates, as predicted from classical theory; i.e.~computed using
  the inverse of the Fisher information. (c) Functional dependence of
  the rescaling constant $\kappa \sigmas^2 / \lambdas$ on the
  parameters $\kappa$ and $\gamma$.}\label{fig:alphasigma}
\end{figure}
\subsection{The average behavior of the MLE}

Our first main result characterizes the `average' behavior of the MLE.
\begin{theorem}\label{thm:mle}
  Assume the dimensionality and signal strength parameters $\kappa$
  and $\gamma$ are such that $\gamma < g_{\text{MLE}}(\kappa)$ (the
  region where the MLE exists asymptotically and shown in Figure
  \ref{fig:kappagamma}). Assume the logistic model described above
  where the empirical distribution of $\{\beta_j\}$ converges weakly
  to a distribution $\Pi$ with finite second moment. Suppose further
  that the second moment converges in the sense that as
  $n \rightarrow \infty$,
  $\operatorname{Ave}_j(\beta_j^{2}) \, \rightarrow \, \E \beta^{2}$,
  $\beta \sim \Pi$.  Then for any pseudo-Lipschitz function $\psi$ of
  order $2$,\footnote{A function $\psi: \R^m \rightarrow \R$ is said
    to be pseudo-Lipschitz of order $k $ if there exists a constant
    $L >0$ such that for all $\bt_0, \bt_1 \in \R^m$,
    $\|\psi(\bt_0) -\psi(\bt_1)\| \leq L\left(1+ \|\bt_0\|^{k-1} +
      \|\bt_1\|^{k-1} \right)\|\bt_0-\bt_1\|$.
  } the marginal distributions of the MLE coordinates obey
\begin{equation}\label{eq:mlenonnull}
  \frac{1}{p} \sum_{j=1}^p \psi(\hbeta_j - \alphas \beta_{j},\beta_j) 
  \,\, {\stackrel{\text{\em a.s.}}{\longrightarrow}} \, \,
  \E [\psi(\sigma_\star Z,\beta)], \quad Z \sim \dnorm (0,1),
\end{equation}
where $\beta \sim \Pi$, independent of $Z$. 
\end{theorem}
\noindent Among the many consequences of this result, we give three:  
\begin{itemize}
\item This result quantifies the exact bias of the MLE in some
  statistical sense. This can be seen by taking $\psi(t,u)= t$ in
  \eqref{eq:mlenonnull}, which leads to
  \[
    \frac{1}{p} \sum_{j=1}^p  (\hbeta_j - \alphas \beta_{j}) \,\,
  {\stackrel{\text{a.s.}}{\longrightarrow}} \, \,
  0, 
\]
and says that $\hat{\beta}_j$ is centered about $\alphas \,
\beta_j$. This can be seen from the empirical results from the
previous sections as well. When $\kappa = 0.2$ and
$\gamma = \sqrt{5}$, the solution to \eqref{eq:main} obeys
$\alphas = 1.499$ and Figure \ref{fig:centering}(a) shows that this is
the correct centering.
\item Second, our result also provides the asymptotic variance of the
  MLE marginals after they are properly centered. This can be seen by
  taking $\psi(t,u) = t^2$, which leads to
  \[
    \frac{1}{p} \sum_{j=1}^p  (\hbeta_j - \alphas \beta_{j})^2 \,\,
  {\stackrel{\text{a.s.}}{\longrightarrow}} \, \,
  \sigmas^2.  
\]
As before, this can also be seen from the empirical results from the
previous section. When $\kappa = 0.2$ and
$\gamma = \sqrt{5}$, the solution to \eqref{eq:main} obeys
$\sigmas = 4.744$ and this is what we see in Figure \ref{fig:varandpval}.
  
\item Third, our result establishes that upon centering the MLE around
  $\alphas \bbeta$, it becomes decorrelated from the signal
  $\bbeta$. This can be seen by taking $\psi(t,u) = t u$, which leads
  to
 \[
    \frac{1}{p} \sum_{j=1}^p  (\hbeta_j - \alphas \beta_{j}) \, \beta_j  \,\,
  {\stackrel{\text{a.s.}}{\longrightarrow}} \, \, 0. 
\]
This can be seen from our earlier empirical results in Figure
\ref{fig:kappagamma}(b). The scatter directly shows the decorrelated
structure and the x-axis passes right through the center,
corroborating our theoretical finding. 
\end{itemize}



It is of course interesting to study how the bias $\alphas$ and the
standard deviation $\sigmas$ depend on the dimensionality $\kappa$ and
the signal strength $\gamma$. We numerically observe that the larger
the dimensionality and/or the larger the signal strength, the larger
the bias $\alphas$. This dependence is illustrated in Figure
\ref{fig:alphasigma}(a). {Further, note that as $\kappa$ approaches
  zero, the bias $\alphas \rightarrow 1$, indicating that the MLE is
  asymptotically unbiased if $p=o(n)$. }The same behavior applies to
$\sigmas$; that is, $\sigmas$ increases in either $\kappa$ or $\gamma$
as shown in Figure \ref{fig:alphasigma}(b). This plot shows the
  theoretical prediction $\sigmas$ divided by the average classical
  standard deviation obtained from $\bm{I}^{-1}({\bbeta})$, the
  inverse of the Fisher information. As $\kappa$ approaches
  zero, the ratio goes to $1$, indicating that the classical standard
  deviation value is valid for $p=o(n)$; this is true across all
  values of $\gamma$. As $\kappa$ increases, the ratio deviates
  increasingly from $1$ and we observe higher and higher variance
  inflation.  In summary, the MLE increasingly deviates from what is
classically expected as either the dimensionality or the signal
strength, or both, increase.

Theorem \ref{thm:mle} is an asymptotic result, and we study how fast
the asymptotic kicks in as we increase the sample size $n$. To this
end, we set $\kappa=0.1$ and let a half of the coordinates of $\bbeta$
have constant value $10$, and the other half be zero. Note that in
this example, $\gamma^2 = 5$ as before. Our goal is to empirically
determine the parameters $\alphas$ and $\sigmas$ from $68,000$ runs,
for each $n$ taking values in $\{2000, 4000, 8000\}$.  Note that there
are several ways of determining $\alphas$ empirically. For instance,
the limit \eqref{eq:mlenonnull} directly suggests taking the ratio
$\sum_j \hat{\beta}_j/\sum_j \beta_j$. An alternative is to consider
taking the ratio when restricting the summation to nonzero indices.
Empirically, we find there is not much difference between these two
choices and choose the latter option, denoting it as $\hat{\alpha}$.
With $\kappa =0.1, \gamma=\sqrt{5}$, the solution to \eqref{eq:main}
is equal to $\alphas = 1.1678, \sigmas = 3.3466,\lambdas= 0.9605$.
Table \ref{tab:alphasigma} shows that $\hat{\alpha}$ is very slightly
larger than $\alphas$ in finite samples.  However, observe that as the
sample size increases, $\hat{\alpha}$ approaches $\alphas$, confirming
the result from \eqref{eq:mlenonnull}. We defer the study of the
asymptotic variance to the next section.
\begin{table}[h!]
  \centering
  \begin{tabular} {|c|c|c|c|}
    \hline
       Parameter & $p=200$ & $p=400$ &$p=800$\\
       \hline
     $\alphas=1.1678$ & $1.1703 (0.0002)$ & $1.1687 (0.0002)$& $1.1681(0.0001)$ \\
     \hline
     $\sigmas = 3.3466$ & $3.3567 (0.0011)$ & $3.3519(0.0008)$ &  $3.3489(0.0006)$ \\
     \hline
    \end{tabular} \\ 
    \caption{Empirical estimates of the centering and standard
      deviation of the MLE.  Standard errors of these estimates are
      between parentheses. In this setting, $\kappa = 0.1$ and
      $\gamma^2 = 5$.  Half of the $\beta_j$'s are equal to ten and
      the others to zero.}
    \label{tab:alphasigma}   
  \end{table}

\subsection{The distribution of the null MLE coordinates}

Whereas Theorem \ref{thm:mle} describes the average or bulk behavior
of the MLE across all of its entries, our next result provides the
explicit distribution of $\hat{\beta}_j$ whenever $\beta_j = 0$,
i.e.~whenever the $j$-th variable is independent from the response
$y$.
\begin{theorem}
  \label{thm:mainthm}
  Let $j$ be any variable such that $\beta_j = 0$. Then in the setting of Theorem \ref{thm:mle}, the MLE obeys  
\begin{equation}\label{eq:firstcoord}
\hbeta_j \, \, {\stackrel{\text{\em d}}{\longrightarrow}} \, \, \dnorm (0,\sigmas^2). 
\end{equation}
For any finite subset of null variables $\{i_1,\hdots, i_k\}$, the
components of $(\hat{\beta}_{i_1},\hdots,\hat{\beta}_{i_k})$ are
asymptotically independent.
\end{theorem}
\noindent In words, the null MLE coordinates are asymptotically normal
with mean zero and variance given by the solution to the system
\eqref{eq:main}. An important remark is this: we have observed that
$\sigmas$ is an increasing function of $\gamma$. Hence, we conclude
that for a null variable $j$, the variance of $\hat{\beta}_j$ is
increasingly larger as the magnitude of the other regression
coefficients increases.

We return to the finite sample precision of the theoretically
predicted asymptotic variance $\sigmas$. As an empirical estimate, we
use $\hat{\beta}_j^2$ averaged over the null coordinates
$\{j : \beta_j = 0\}$ since it is approximately unbiased for
$\sigmas^2$. We work in the setting of Table \ref{tab:alphasigma} in
which $\sigmas = 3.3466$, averaging our $68,000$ estimates.  The
results are given in this same table; we observe that $\hat{\sigma}$
is very slightly larger than $\sigmas$. However, it progressively gets
closer to $\sigmas$ as the sample size $n$ increases.

Next, we study the accuracy of the asymptotic convergence results in
\eqref{eq:firstcoord}.  In the setting of Table \ref{tab:alphasigma},
we fit $500,000$ independent logistic regression models and plot the
empirical cumulative distribution function of
$\Phi(\hbbeta_j/\sigmas)$ in Figure \ref{fig:pvaladjusted}(a) for some
fixed null coordinate. Observe the perfect agreement with a straight
line of slope 1.

\subsection{The distribution of the LRT}

We finally turn our attention to the distribution of the likelihood
ratio statistic for testing $\beta_j = 0$.
\begin{theorem}\label{thm:lrt}
  Consider the LLR
  $\Lambda_j = \min_{\bb \, : \, b_j = 0} \ell(\bb) - \min_{\bb}
  \ell(\bb)$ for testing $\beta_j = 0$. In the setting of Theorem
  \ref{thm:mle}, twice the LLR is asymptotically distributed as a
  multiple of a chi-square under the null,
\begin{equation}\label{eq:lrt}
2 \Lambda_j   \, \, {\stackrel{\text{\em d}}{\longrightarrow}} \, \,  \frac{\kappa \, \sigmas^2}{\lambdas} \, \chi^2_{1}. 
\end{equation}
Also, the LLR for testing $\beta_{i_1}=\beta_{i_2} =\hdots= \beta_{i_k} =
  0$ for any finite $k$ converges to the rescaled
  chi-square $\left(\kappa \sigmas^2/\lambdas\right) \chi^2_k$ under the null.
\end{theorem}
This theorem explicitly states that the LRT does not follow a
$\chi_1^2$ distribution as soon as $\kappa > 0$ since the
multiplicative factor is then larger than one, as demonstrated in
Figure \ref{fig:alphasigma}(c). In other words, the LRT is
stochastically quite larger than a $\chi_1^2$, explaining the large
spike near zero in Figure \ref{fig:LRT}. Also, Figure
\ref{fig:alphasigma}(c) suggests that as $\kappa \rightarrow 0$, the
classical result is recovered.

Theorem \ref{thm:lrt} extends to arbitrary signal strengths the
earlier result from \cite{sur2017likelihood}, which described the
distribution of the LRT under the global null $(\beta_j = 0$ for all
$j$). One can quickly verify that when $\gamma = 0$, the
multiplicative factor in \eqref{eq:lrt} is that given in
\cite{sur2017likelihood}, which easily follows from the fact that in
this case, the system \eqref{eq:main} reduces to
\eqref{eq:reduced}. {Furthermore, if the signal is sparse in the sense
  that $o(n)$ coefficients have non-zero values, $\gamma^2=0$, which
  immediately implies that the asymptotic distribution for the LLR
  from \cite{sur2017likelihood} still holds in such cases.}

\begin{figure}
        \centering
        \begin{subfigure}[b]{0.35\textwidth}
           \centering
           \hspace{-0.5cm}
            \includegraphics[width=\textwidth]{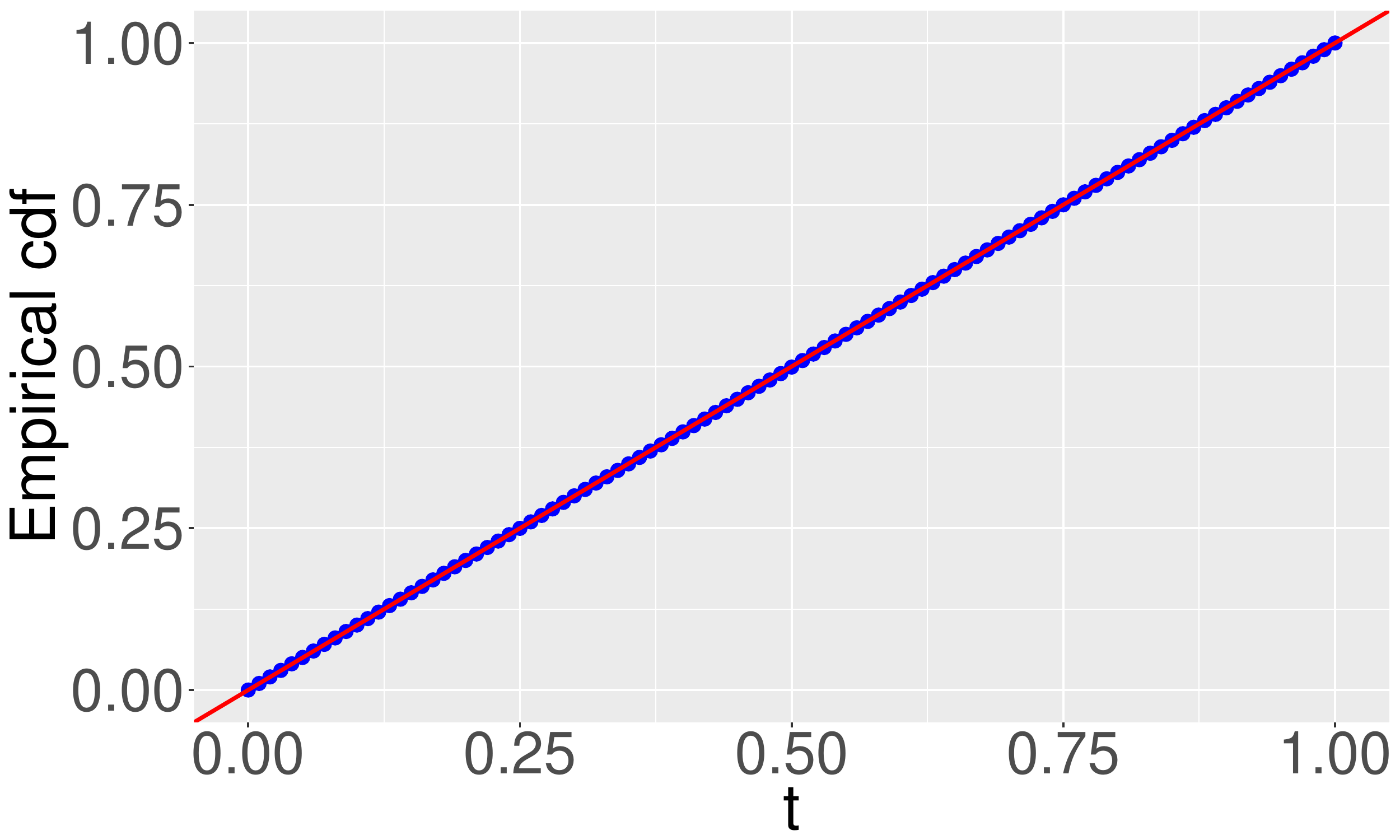}
           \caption{}
        \end{subfigure}
        \begin{subfigure}[b]{0.35\textwidth}  
            \centering 
            \includegraphics[width=\textwidth]{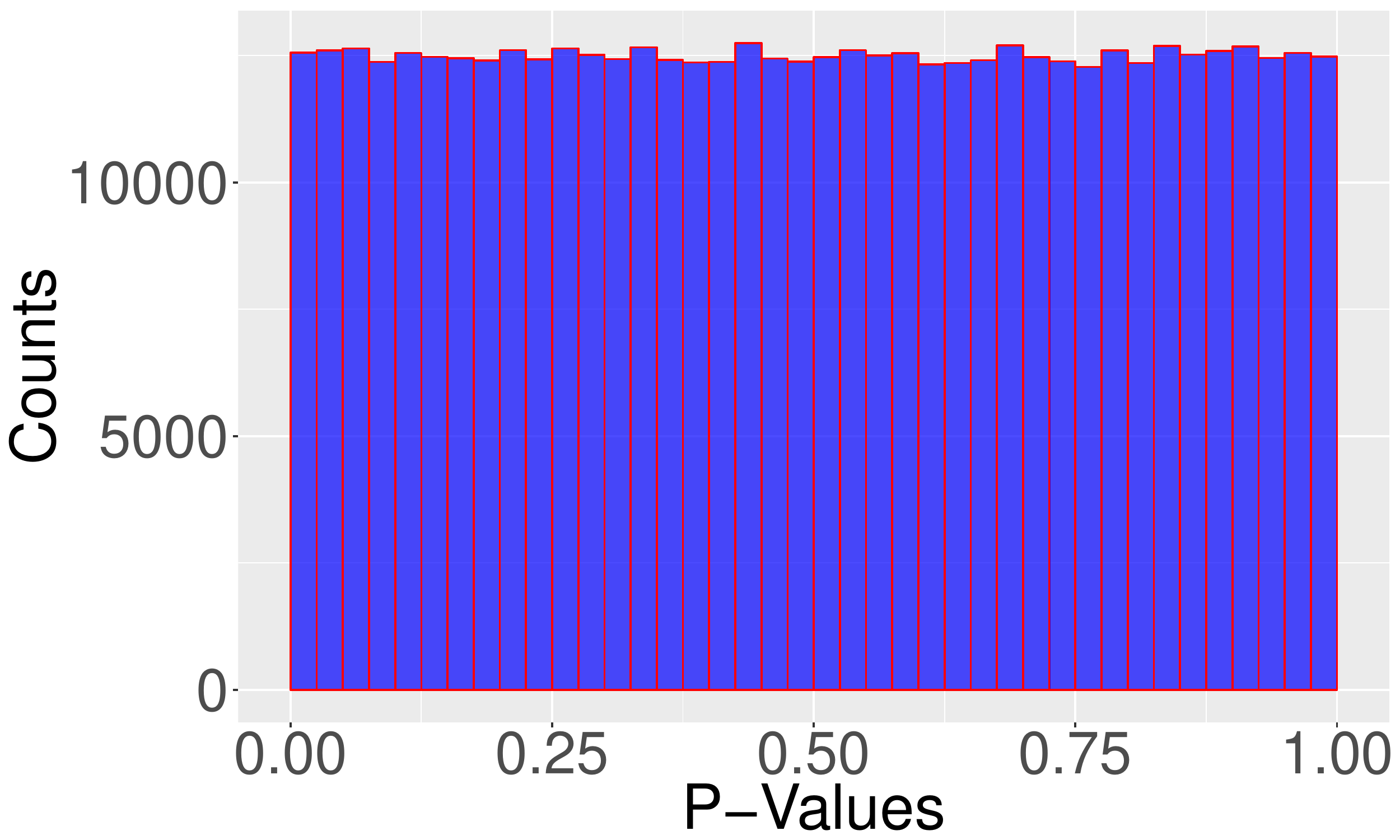}
            \caption{}    
            
        \end{subfigure}
        \vskip\baselineskip
        \begin{subfigure}[b]{0.35\textwidth}   
            \centering
            \hspace{-0.5cm} 
            \includegraphics[width=\textwidth]{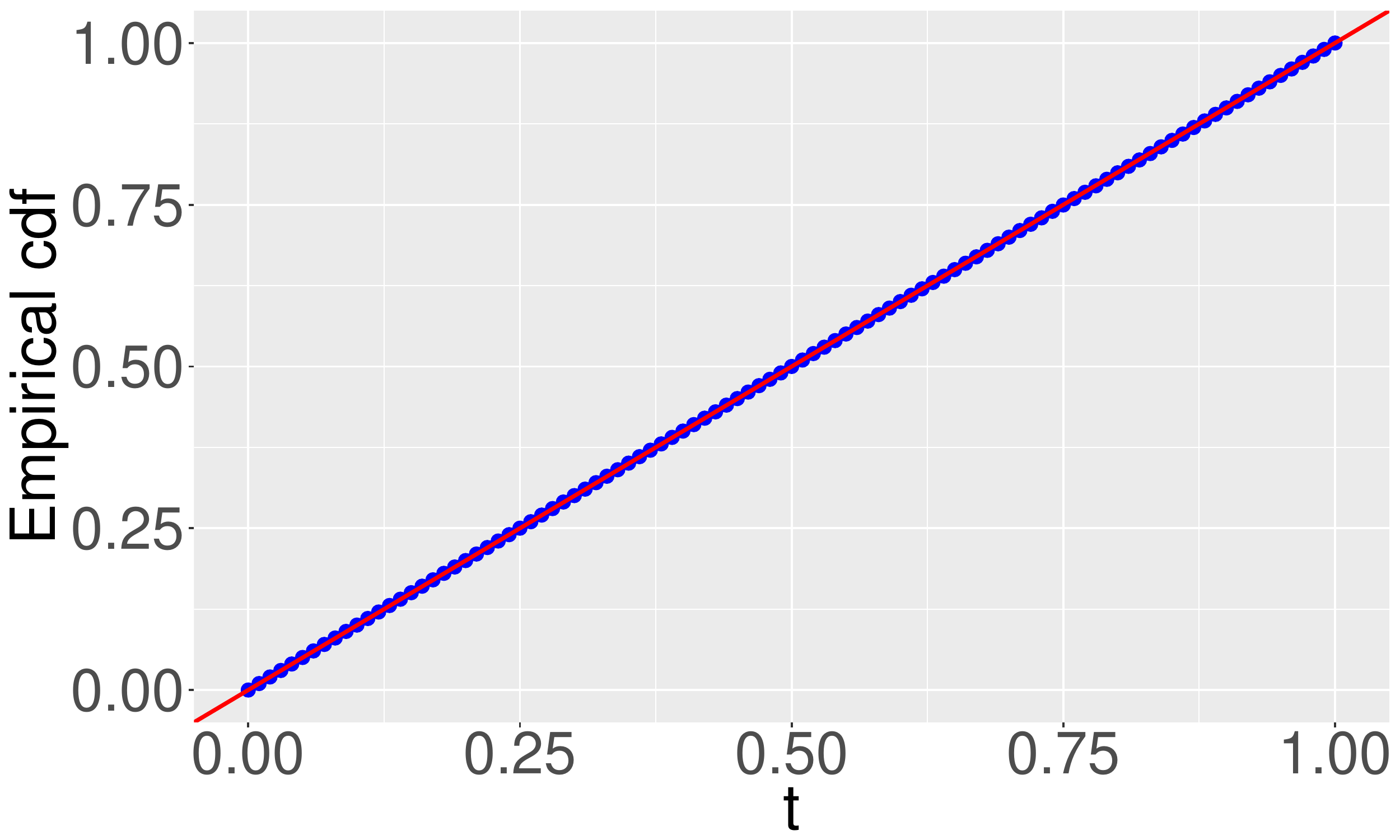}
            \caption{}    
               \end{subfigure}
        \begin{subfigure}[b]{0.35\textwidth}   
            \centering 
            \includegraphics[width=\textwidth]{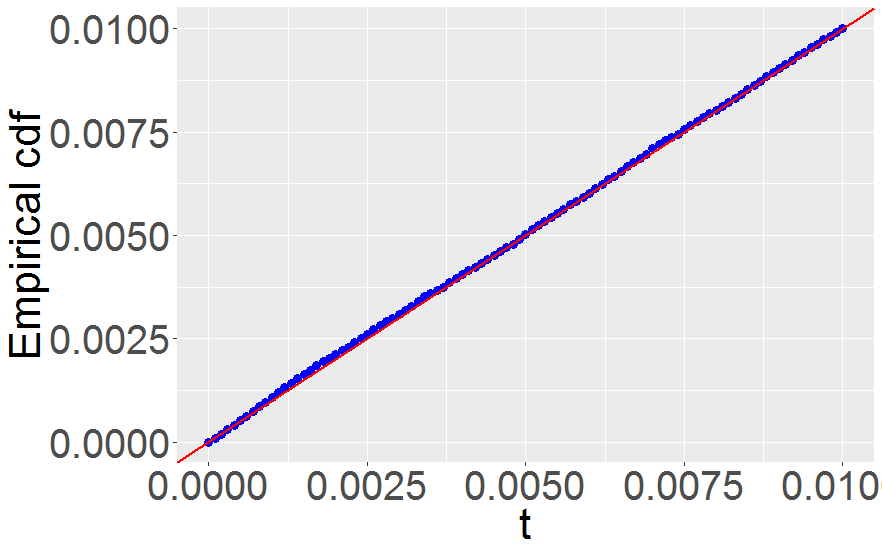}
            \caption{}    
        \end{subfigure}
        \caption{The setting is that from Table \ref{tab:alphasigma} with
  $n = 4000$.  (a) Empirical cdf of $\Phi(\hat{\beta}_j/\sigmas)$ for
  a null variable ($\beta_j = 0$).  (b) P-values given by the LLR
  approximation \eqref{eq:lrt} for this same null variable.  (c)
  Empirical distribution of the p-values from (b). (d) Same as (c) but
  showing accuracy in the lower tail (check the range of the
  horizontal axis). All these plots are based on 500,000
  replicates.}
       \label{fig:pvaladjusted}
    \end{figure}

To investigate the quality of the accuracy of \eqref{eq:lrt} in finite
samples, we work on the p-value scale. We select a null coefficient
and compute p-values based on \eqref{eq:lrt}.  The histogram for the
p-values across $500,000$ runs is shown in Figure
\ref{fig:pvaladjusted}(b) and the empirical cumulative distribution
function (cdf) in Figure \ref{fig:pvaladjusted}(c).  In stark contrast to
Figure \ref{fig:varandpval}, we observe that the p-values are uniform
over the bulk of the distribution.

From a multiple testing perspective, it is essential to understand the
accuracy of the rescaled chi-square approximation in the tails of the
distribution. We plot the empirical cdf of the p-values, zooming in
the tail, in Figure \ref{fig:pvaladjusted}(d). We find that the
rescaled chi-squared approximation works extremely well even in the
tails of the distribution.

\begin{table}[h!]
\centering
    \begin{tabular} {|c|c|c|}
    \hline
     & $p=400$ & $p=800$ \\ \hline
    $\prob\{\text{p-value}  \leq 5 \% \}$ & $5.03 \% (0.031 \%)$ & $5.01 \% (0.03 \%)$ \\ 
    \hline
        $\prob\{\text{p-value} \leq 1 \% \}$ & $1.002\% (0.014 \%)$ & $1.005 \% (0.014 \%)$  \\ 
        \hline
        $\prob\{\text{p-value}  \leq 0.5 \% \}$& $0.503 \% (0.01 \%)$ & $0.49\% (0.0099 \%)$ \\ 
        \hline
        $\prob\{\text{p-value}  \leq 0.1 \% \}$ &  $0.109 \% (0.004 \%)$ & $0.096 \% (0.0044 \%)$\\ \hline
        $\prob\{\text{p-value}  \leq 0.05 \% \}$ &  $0.052 \% (0.003 \%)$ & $0.047 \% (0.0031 \%)$ \\ \hline
        $\prob\{\text{p-value}  \leq 0.01 \% \}$ &  $0.008 \% (0.0013 \%)$ & $0.008 \% (0.0013 \%)$ \\ \hline
    \end{tabular} \\ 
    \caption{P-value probabilities estimated over $500,000$
      replicates with standard errors in parentheses. Here, $\kappa =
      0.1$ and the  setting is otherwise the same as in Table \ref{tab:alphasigma}.}
    \label{tab:pvaladjustedfits}   
\end{table}
To obtain a more refined idea of the quality of approximation, we zoom
in the smaller bins close to zero and provide estimates of the p-value
probabilities in Table \ref{tab:pvaladjustedfits} for $n=4000$ and
$n=8000$. The tail approximation is accurate, modulo a slight
deviation in the bin for $\prob\{\text{p-value}\}  \leq 0.1$ for the smaller sample size. For $n=8000$,
however, this deviation vanishes and we find perfect coverage of the
true values. It seems that our approximation is extremely precise even
in the tails.

\subsection{Other scalings}

Throughout this section, we worked under the assumption that
  $\lim_{n \rightarrow \infty} \text{Var}(\bX_i'\bbeta) = \gamma^2$,
  which does not depend on $n$, and we explained that this is the only
  scaling that makes sense to avoid a trivial problem. We set the
  variables to have variance $1/n$ but this is of course somewhat
  arbitrary. For example, we could choose them to have variance $v$ as
  in $\bX_i \sim \dnorm(\bzero, v\bm{I}_{p})$. This means that
  $\bX_i = \sqrt{vn} \bZ_i$, where $\bZ_i$ is as before. This gives
  $ \bX_i' \bbeta = \bZ_i' \bb$, where $\bbeta = \bb/\sqrt{n v}$.  The
  conclusions from Theorem \ref{thm:mle} and \ref{thm:mainthm} then
  hold for the model with predictors $\bZ_i$ and regression
  coefficient sequence $\bb$. Consequently, by simple rescaling, we
  can pass the properties of the MLE in this model to those of the MLE
  in the model with predictors $\bX_i$ and coefficients $\bbeta$. For
  instance, the asymptotic standard deviation of $\hbbeta$ is equal to
  $\sigmas/\sqrt{nv}$, where $\sigmas$ is just as in Theorems
  \ref{thm:mle} and \ref{thm:mainthm}. On the other hand, the result
  for the LRT, namely, Theorem \ref{thm:lrt} is scale invariant; no
  such trivial adjustment is necessary.

\subsection{Non-Gaussian covariates}
\label{sec:nonGaussian}
Our model assumes that the features are Gaussian, yet, we expect that
the same results hold under other distributions with the proviso that
they have sufficiently light tails.  In this section, we empirically
study the applicability of our results for certain non-Gaussian
features.

\begin{figure}
\begin{center}
	\begin{tabular}{cc}
	\includegraphics[scale=0.3,keepaspectratio]{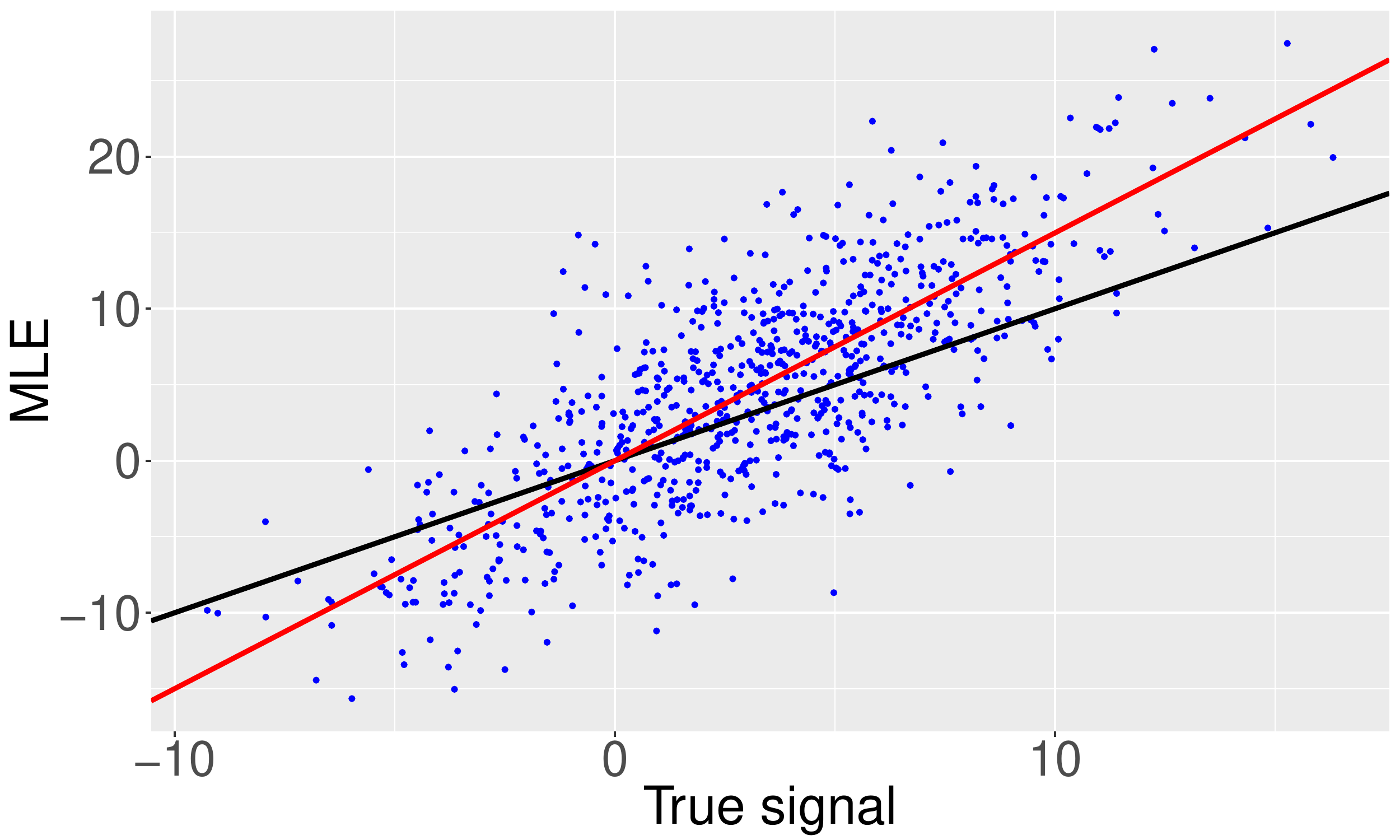} 
		&  \includegraphics[scale=0.3,keepaspectratio]{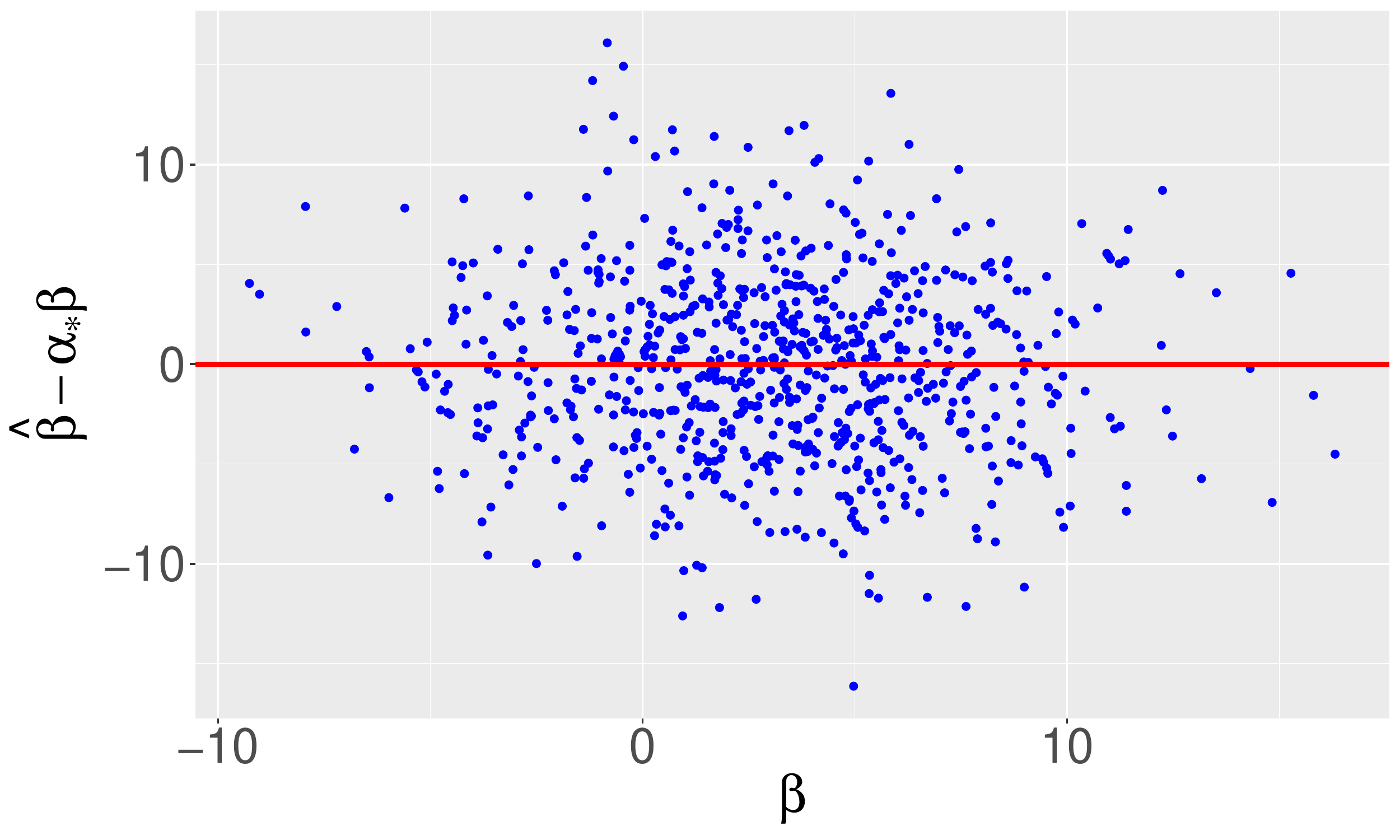}\label{fig:var} \tabularnewline
		(a) & (b)  \tabularnewline
	\end{tabular}
\end{center}
\caption{Simulation for a non-Gaussian design. The $j$-th feature
  takes values in $\{0,1,2\}$ with probabilities
  $p_j^2, 2p_j(1-p_j), (1-p_j)^2$; here, $p_j \in [0.25,0.75]$ and
  $p_j \neq p_k$ for $j \neq k$. Features are then centered and
  rescaled to have unit variance. The setting is otherwise the same as
  for Figure \ref{fig:centering}. (a) Analogue of Figure
  \ref{fig:centering}(a). Red line has slope
  $\alphas \approx 1.499$.(b) Analogue of Figure
  \ref{fig:kappagamma}(b). Observe the same behavior as earlier: the
  theory predicts correctly the bias and the decorrelation between the
  bias-adjusted residuals and the true effect sizes.}
\label{fig:centering_SNPdesign}
\end{figure}
In genetic studies, we often wish to understand how a binary
response/phenotype depends on single nucleotide polymorphisms (SNPs),
which typically take on values in $\{0,1,2\}$. When the $j$-th SNP is
in Hardy-Weinberg equilibrium, the chance of observing $0$, $1$ and
$2$ is respectively $p_j^2$, $2p_j(1-p_j)$ and $(1-p_j)^2$, where
$p_j$ is between $0$ and $1$. Below we generate independent features
with marginal distributions as above for parameters $p_j$ varying in
$[0.25, 0.75]$.  We then center and normalize each column of the
feature matrix $\bX$ to have zero mean and unit variance. Keeping
everything else as in the setting of Figure \ref{fig:centering}, we
study the bias of the MLE in Figure
\ref{fig:centering_SNPdesign}(a). As for Gaussian designs, the MLE
seriously over-estimates effect magnitudes and our theoretical
prediction $\alphas$ accurately corrects for the bias. We also see
that the bias-adjusted residuals $\hbbeta - \alphas \bbeta$ are
uncorrelated with the effect sizes $\bbeta$, as shown in Figure
\ref{fig:centering_SNPdesign}(b).

 \begin{figure}
        \centering
        \begin{subfigure}[b]{0.35\textwidth}
           \centering
           \hspace{-0.5cm}
            \includegraphics[width=\textwidth]{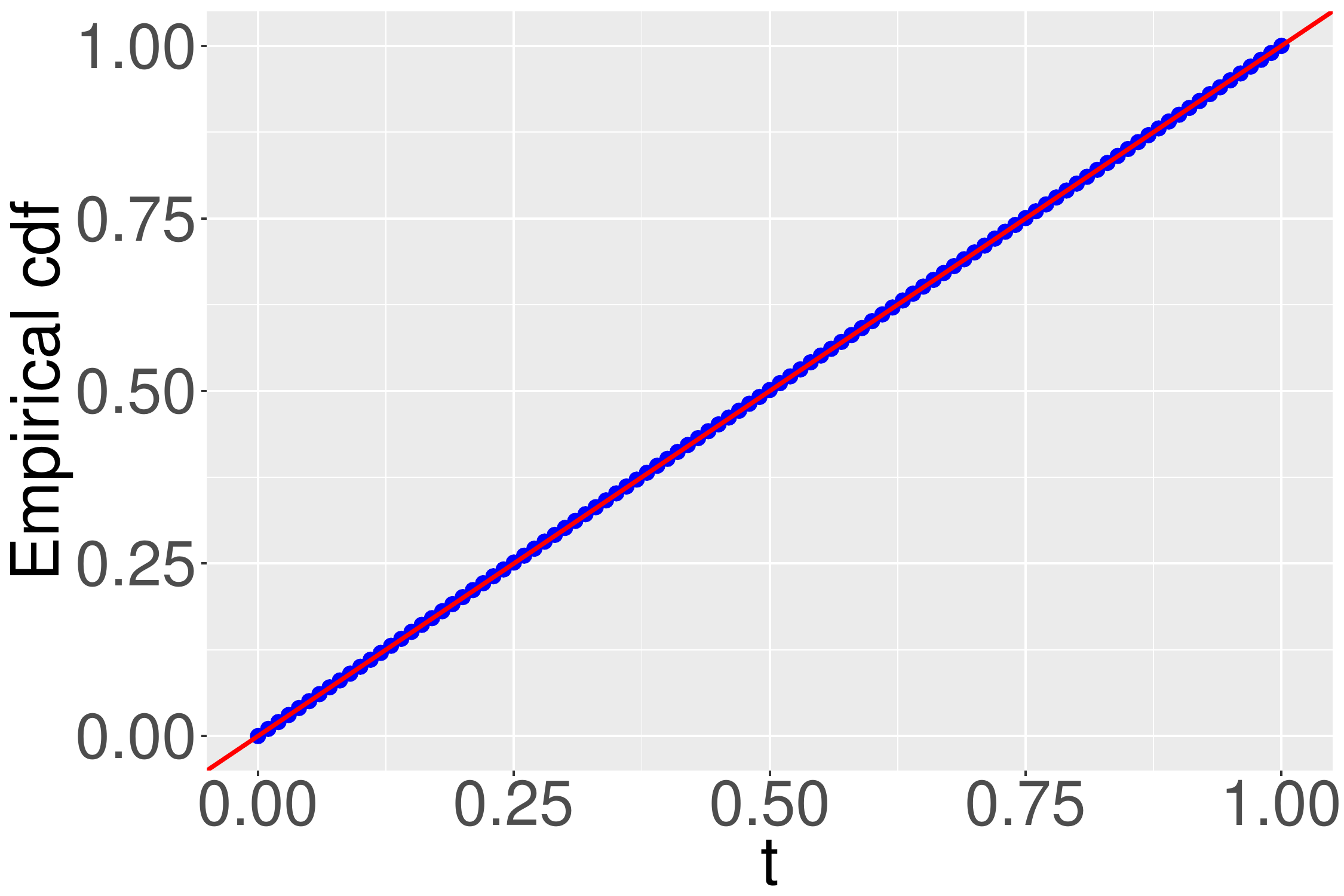}
           \caption{}
        \end{subfigure}
        \begin{subfigure}[b]{0.35\textwidth}  
            \centering 
            \includegraphics[width=\textwidth]{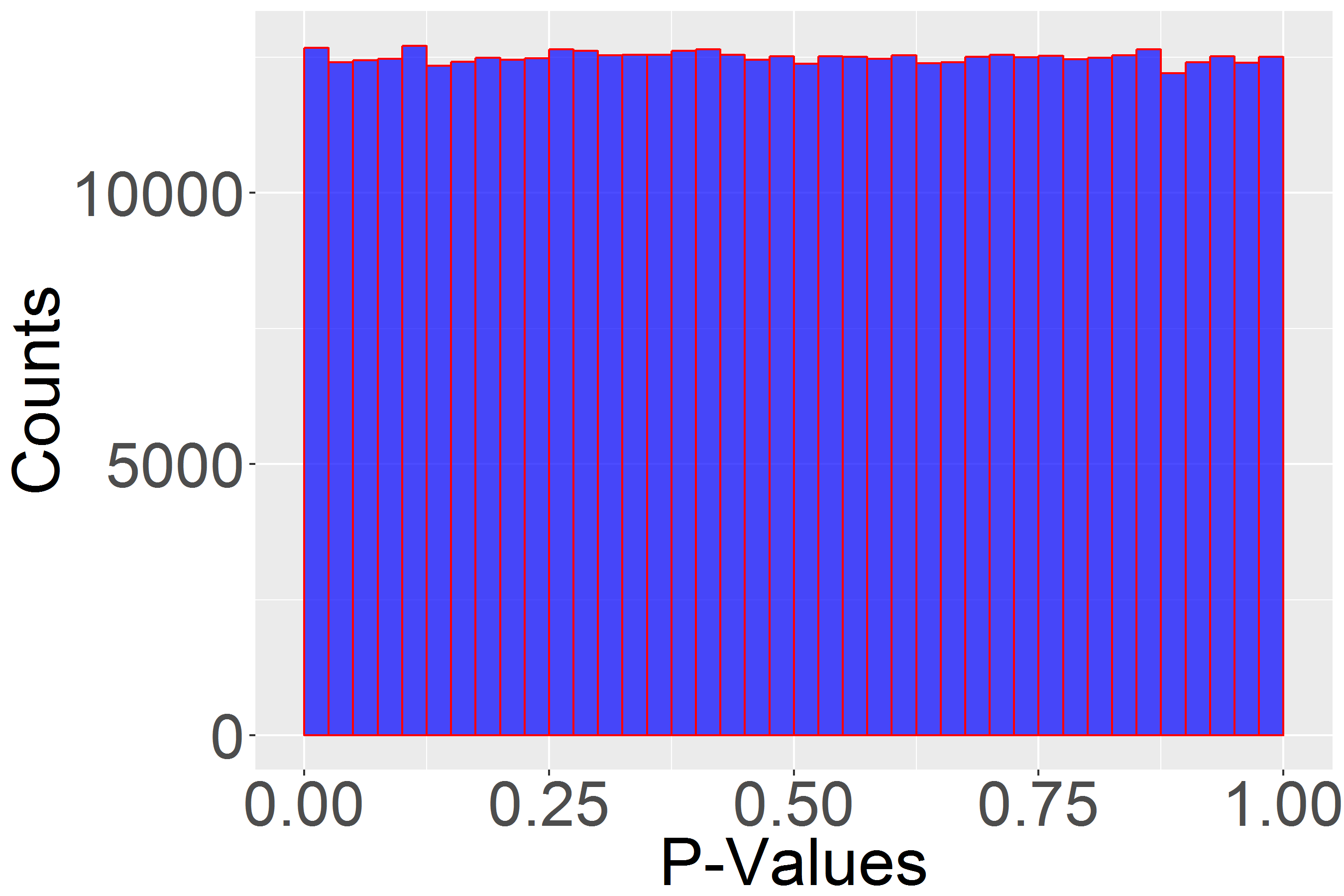}
            \caption{}    
            
        \end{subfigure}
        \vskip\baselineskip
        \begin{subfigure}[b]{0.35\textwidth}   
            \centering
            \hspace{-0.5cm} 
            \includegraphics[width=\textwidth]{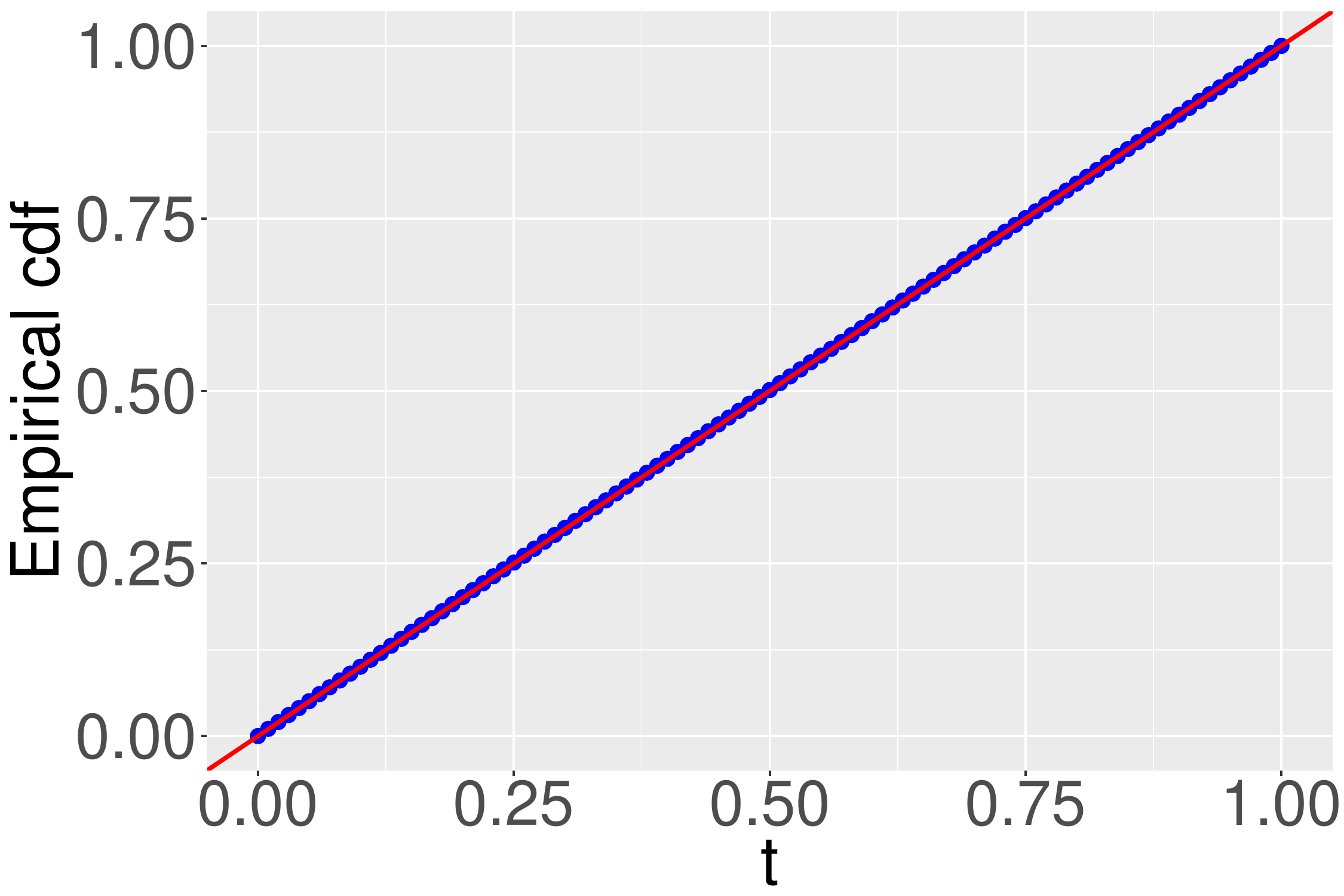}
            \caption{}    
               \end{subfigure}
        \begin{subfigure}[b]{0.35\textwidth}   
            \centering 
            \includegraphics[width=\textwidth]{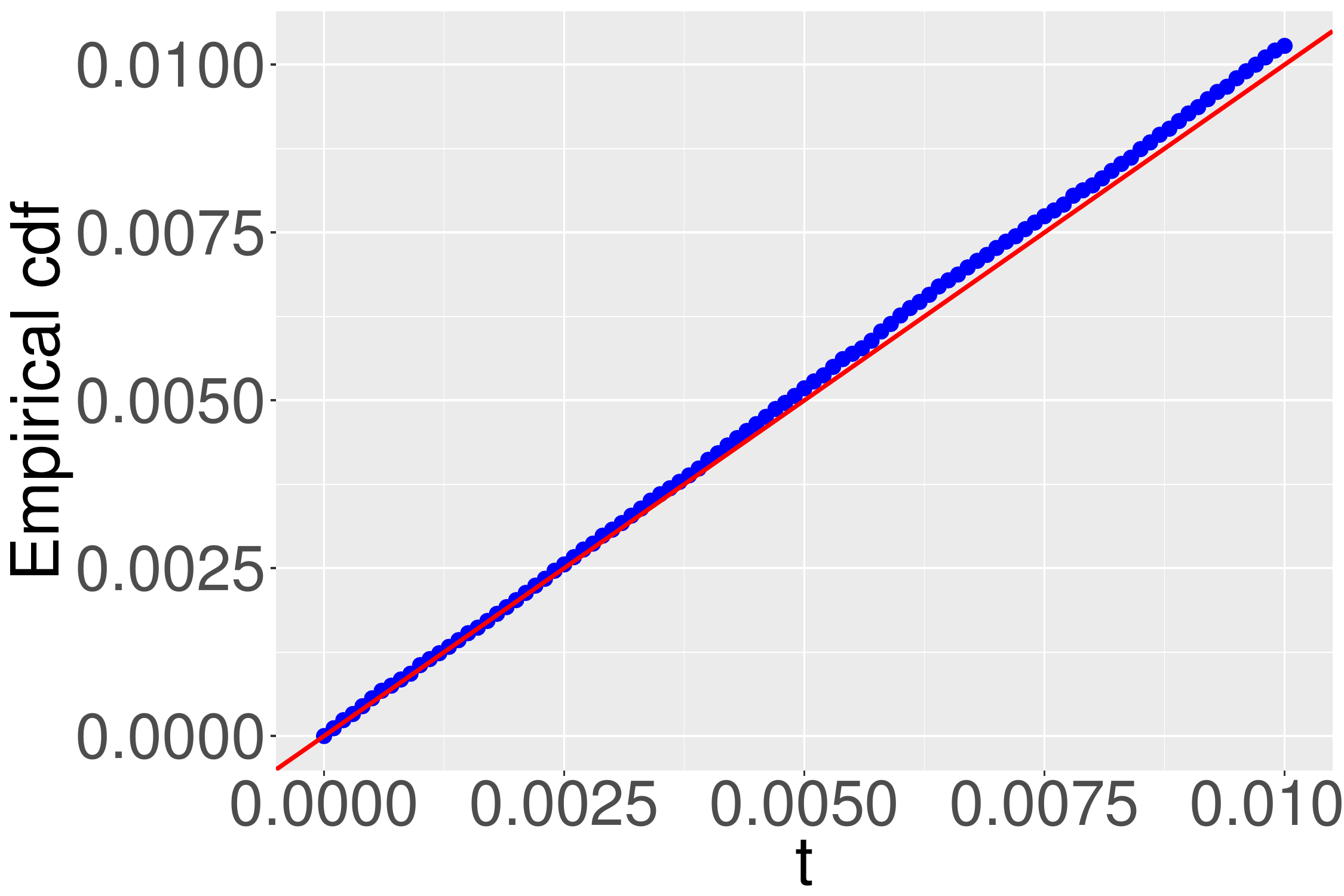}
            \caption{}    
        \end{subfigure}
        \caption{The features are multinomial as in Figure
          \ref{fig:centering_SNPdesign} and the setting is otherwise
          the same as for Figure \ref{fig:pvaladjusted}.(a) Empirical
          cdf of $\Phi(\hat{\beta}_j/\sigmas)$ for a null variable
          ($\beta_j = 0$).  (b) P-values given by the LLR
          approximation \eqref{eq:lrt} for this same null variable.
          (c) Empirical distribution of the p-values from (b). (d)
          Same as (c) but displaying accuracy in the extreme. These
          results are based on $500,000$
          replicates.}  \label{fig:snpdesign}
\end{figure}
The bulk distribution of a null coordinate suggested by Theorem
\ref{thm:mainthm} and the LRT distribution from Theorem \ref{thm:lrt}
are displayed in Figure \ref{fig:snpdesign}. Other than the design,
the setting is the same as for Figure \ref{fig:pvaladjusted}. The
theoretical predictions are once again accurate.  Furthermore, upon
examining the tails of the p-value distribution, we once more observe
a close agreement with our theoretical predictions.  All in all,
these findings indicate that our theory is expected to apply to a far
broader class of features.

\section{Adjusting Inference by Estimating  the Signal Strength}\label{sec:gammaest}

All of our asymptotic results, namely, the average behavior of the
MLE, the asymptotic distribution of a null coordinate, and the LLR,
depend on the unknown signal strength $\gamma$.  In this section, we
describe a simple procedure for estimating this single parameter from
an idea proposed by Boaz Nadler and Rina Barber after the second
author presented the new high-dimensional ML theory from this paper at
the Mathematisches Forshunginstitut Oberwolfach on March 12, 2018.

\newcommand{\gMLE}{g_{\text{MLE}}}

\subsection{ProbeFrontier: estimating $\gamma$ by probing the MLE
  frontier}

We estimate the signal strength by actually using the predictions from
our theory, namely, the fact that we have asymptotically characterized
in Section 2 the region where the MLE exists. We know from Theorem
\ref{thm:mle_exist} that for each $\gamma$, there is a maximum
dimensionality $\gMLE^{-1}(\gamma)$ at which the MLE ceases to
exist. We propose an estimate $\hat \kappa$ of $\gMLE^{-1}(\gamma)$
and set $\hat \gamma = \gMLE(\hat \kappa)$. Below, we shall refer to
this as the {\em ProbeFrontier} method.

 Given a data sample $(y_i,\bX_i)$, we begin by choosing a fine
  grid of values
  $\kappa \le \kappa_1 \le \kappa_2 \le ... \le \kappa_K \le 1/2$. For
  each $\kappa_j$, we execute the following procedure:
 \begin{description}
 \item[{\em Subsample}] Sample $n_j = p/\kappa_j$
   observations from the data without replacement, rounding to the
   nearest integer. Ignoring the rounding, the dimensionality of this
   subsample is $p/n_j = \kappa_j$.
 \item[{\em Check whether MLE exists}] For the subsample, check
   whether the MLE exists or not. This is done by solving a linear
   programming feasibility problem; if there exists a vector
   $\bm{b} \in \R^p$ such that $\bX_i'\bm{b}$ is positive when
   $y_i=1$ and negative otherwise, then perfect separation between
   cases and controls occurs and the MLE does not exist.  Conversely,
   if the linear program is infeasible, then the MLE exists.
 \item[{\em Repeat}] Repeat the two previous steps $B$ times and
   compute the proportion  of times $\hat{\pi}(\kappa_j)$ the MLE does
   not exist.
\end{description}
We next find $(\kappa_{j-1},\kappa_j)$, such that $\kappa_j$ is the
smallest value in $\mathcal{K}$ for which $\hat{\pi}(\kappa_j) \ge 0.5$.
By linear interpolation
between $\kappa_{j-1}$ and $\kappa_j$, we obtain $\hat{\kappa}$ for
which the proportion of times the MLE does not exist would be $0.5$. 
We set $\hat \gamma = \gMLE(\hat \kappa)$.
(Since the
`phase-transition' boundary for the existence of the MLE is a smooth
function of $\kappa$, as is clear from Figure \ref{fig:kappagamma},
choosing a sufficiently fine grid $\{\kappa_j\}$ would make the linear
interpolation step sufficiently precise.)  

\subsection{Empirical performance of adjusted inference}
\begin{figure}
  \centering
  \includegraphics[scale=0.35,keepaspectratio]{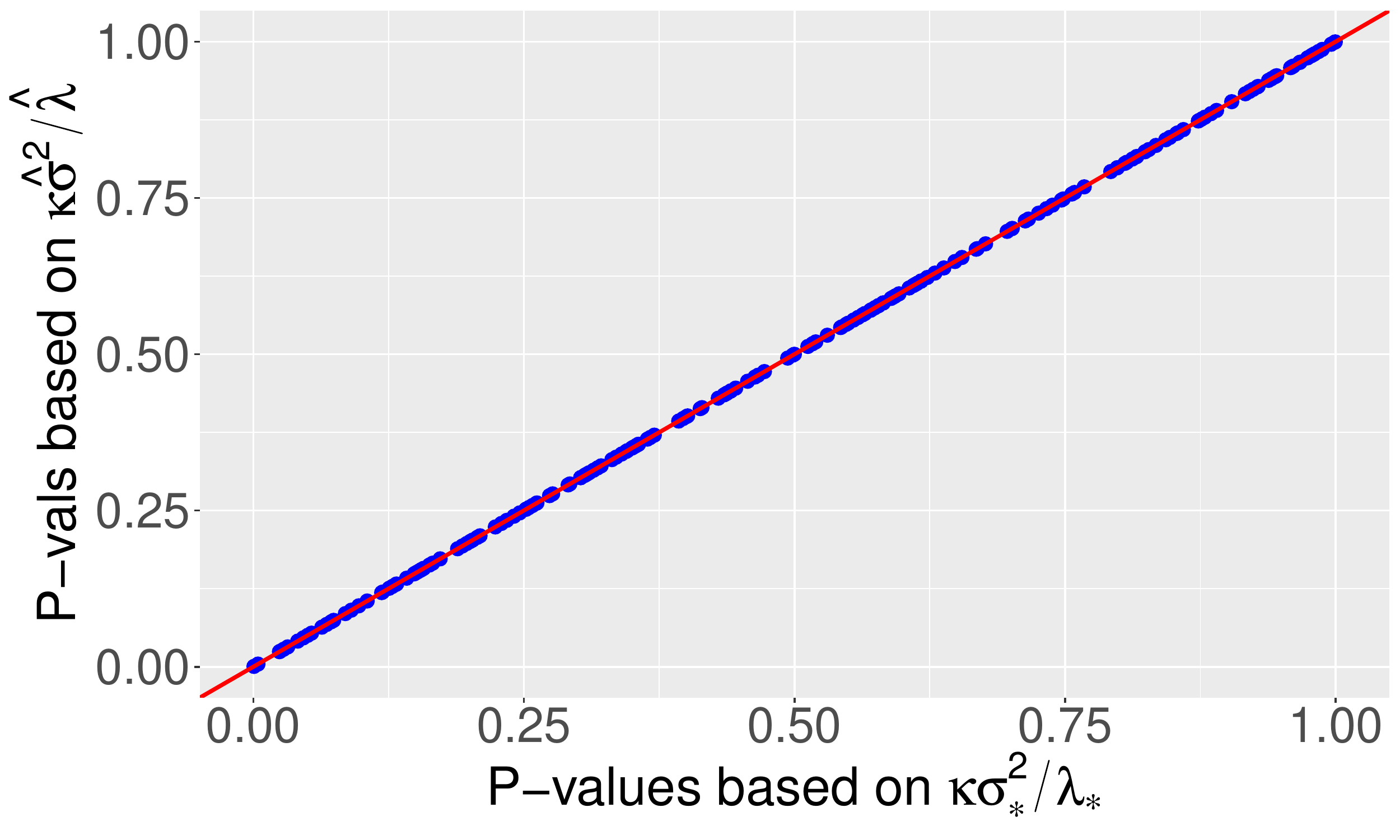}  
  \caption{(a) Null p-values obtained using the
    $(\kappa \hat{\sigma}^2/\hat{\lambda}) \, \chi_1^2$ approximation
    plotted against those obtained using
    $(\kappa \sigmas^2/\lambdas) \, \chi_1^2$. Observe the perfect
    agreement (the diagonal is in red).}
\label{fig:gamma_hat_prec}
\end{figure}  
We demonstrate the accuracy of ProbeFrontier via some empirical
results. We begin by generating $4000$ i.i.d.~observations
$(y_i,\bX_i)$ using the same setup as in Figure \ref{fig:pvaladjusted}
($\kappa=0.1$ and half of the regression coefficients are null). We
work with a sequence $\{\kappa_j\}$ of points spaced apart by
$10^{-3}$ and obtain $\hat{\gamma}$ via the procedure described above,
drawing $50$ subsamples. Solving the system \eqref{eq:main} using
$\kappa=0.1$ and $\hat{\gamma}$ yields estimates for the theoretical
predictions $(\alphas,\sigmas,\lambdas)$ equal to
$(\hat \alpha, \hat \sigma, \hat \lambda)$ $= (1.1681,3.3513,0.9629)$.
In turn, this yields an estimate for the multiplicative factor
$\kappa \sigmas^2/\lambdas$ in \eqref{eq:lrt} equal to $ 1.1663$.
Recall from Section \ref{sec:main} that the theoretical values are
$(\alphas,\taus,\lambdas) = (1.1678,3.3466,0.9605)$ and
$\kappa \sigmas^2/\lambdas = 1.1660$. Next, we compute the LLR
statistic for each null and p-values from the approximation
\eqref{eq:lrt} in two ways: first, by using the theoretically
predicted values, and second, by using our estimates. A scatterplot of
these two sets of p-values is shown in Figure
\ref{fig:gamma_hat_prec}(a) (blue). We observe impeccable agreement.

\begin{figure}
        \centering
        \begin{subfigure}[b]{0.35\textwidth}
           \centering
           \hspace{-0.5cm}
            \includegraphics[width=\textwidth]{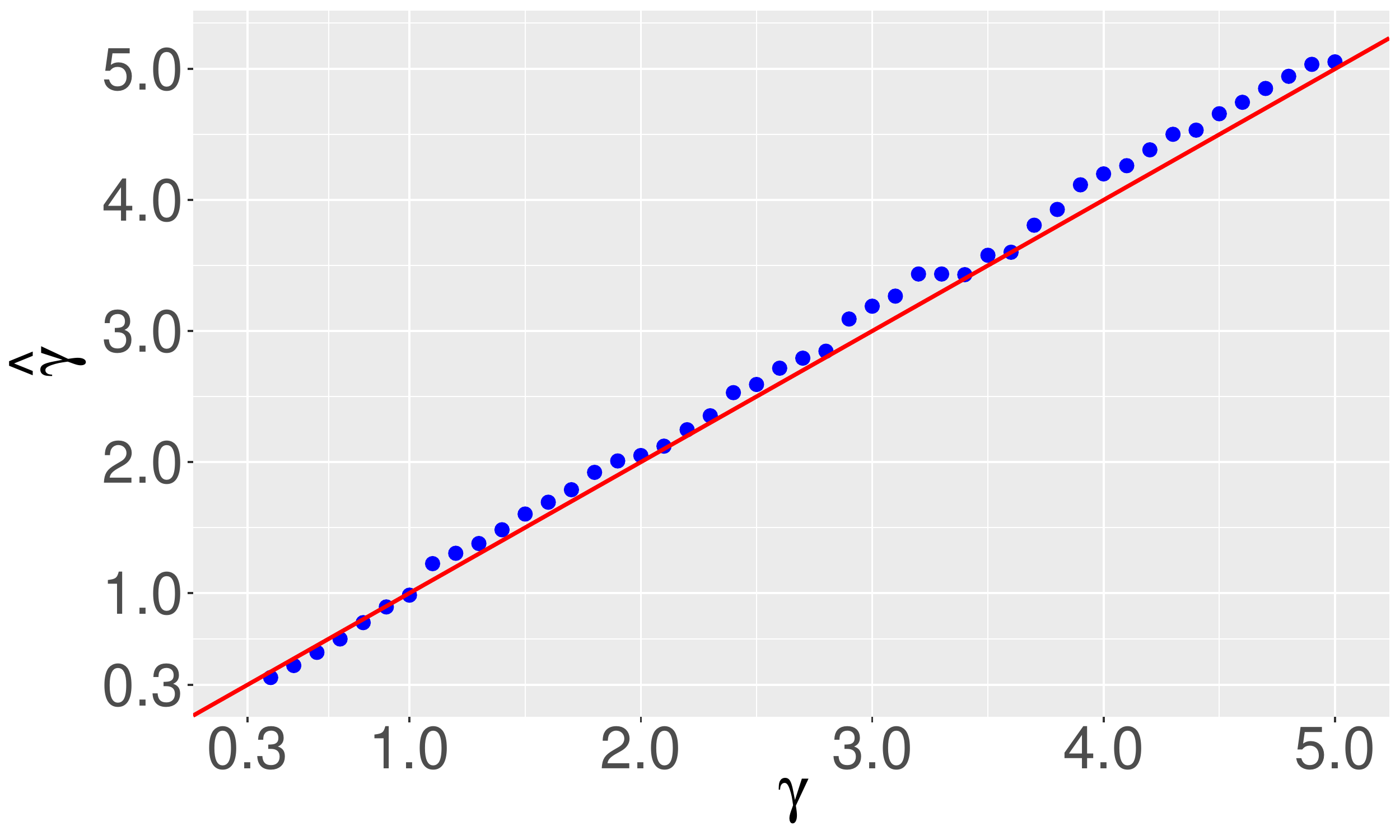}
           \caption{}
        \end{subfigure}
        \begin{subfigure}[b]{0.35\textwidth}  
            \centering 
            \includegraphics[width=\textwidth]{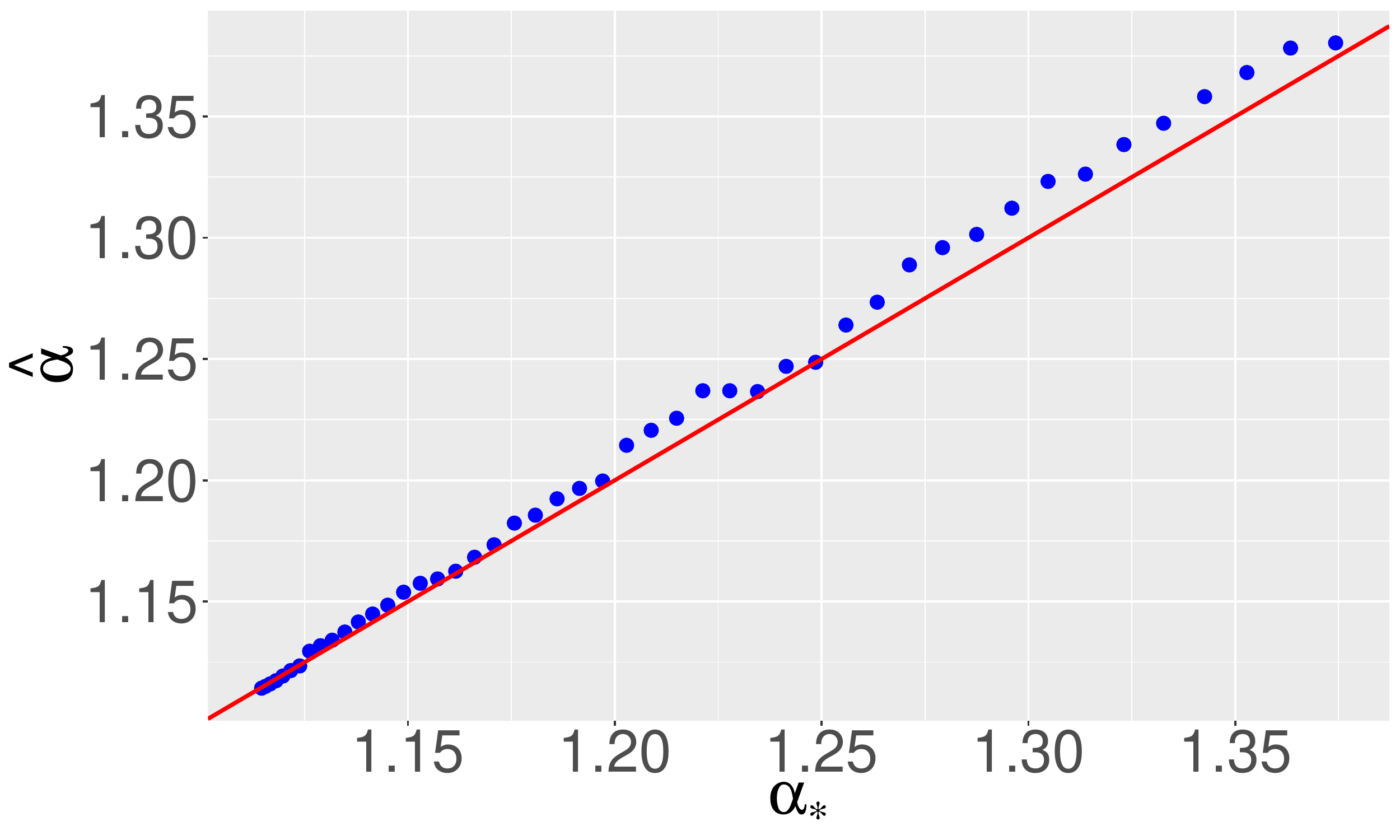}
            \caption{}    
            
        \end{subfigure}
        \vskip\baselineskip
        \begin{subfigure}[b]{0.35\textwidth}   
            \centering
            \hspace{-0.5cm} 
            \includegraphics[width=\textwidth]{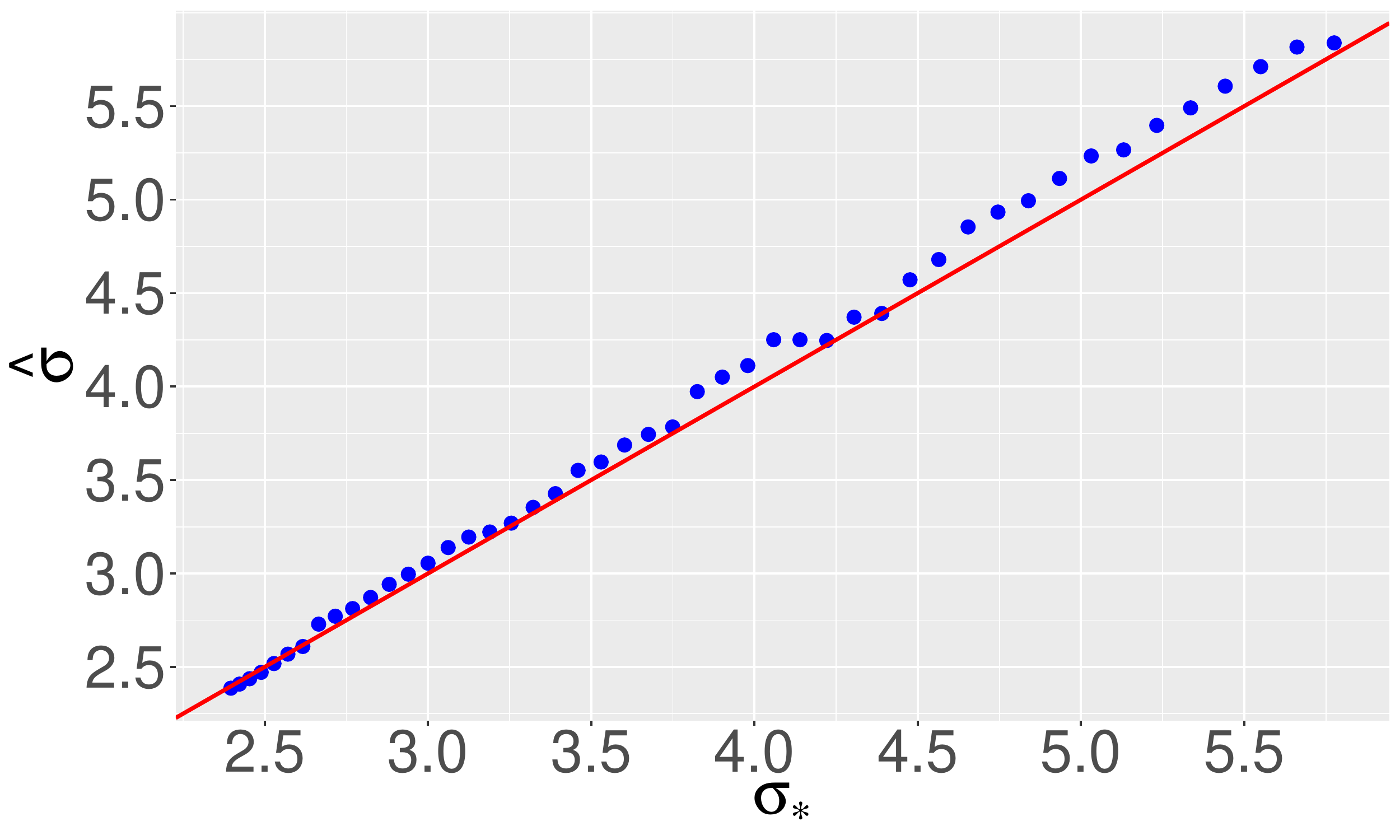}
            \caption{}    
               \end{subfigure}
        \begin{subfigure}[b]{0.35\textwidth}   
            \centering 
            \includegraphics[width=\textwidth]{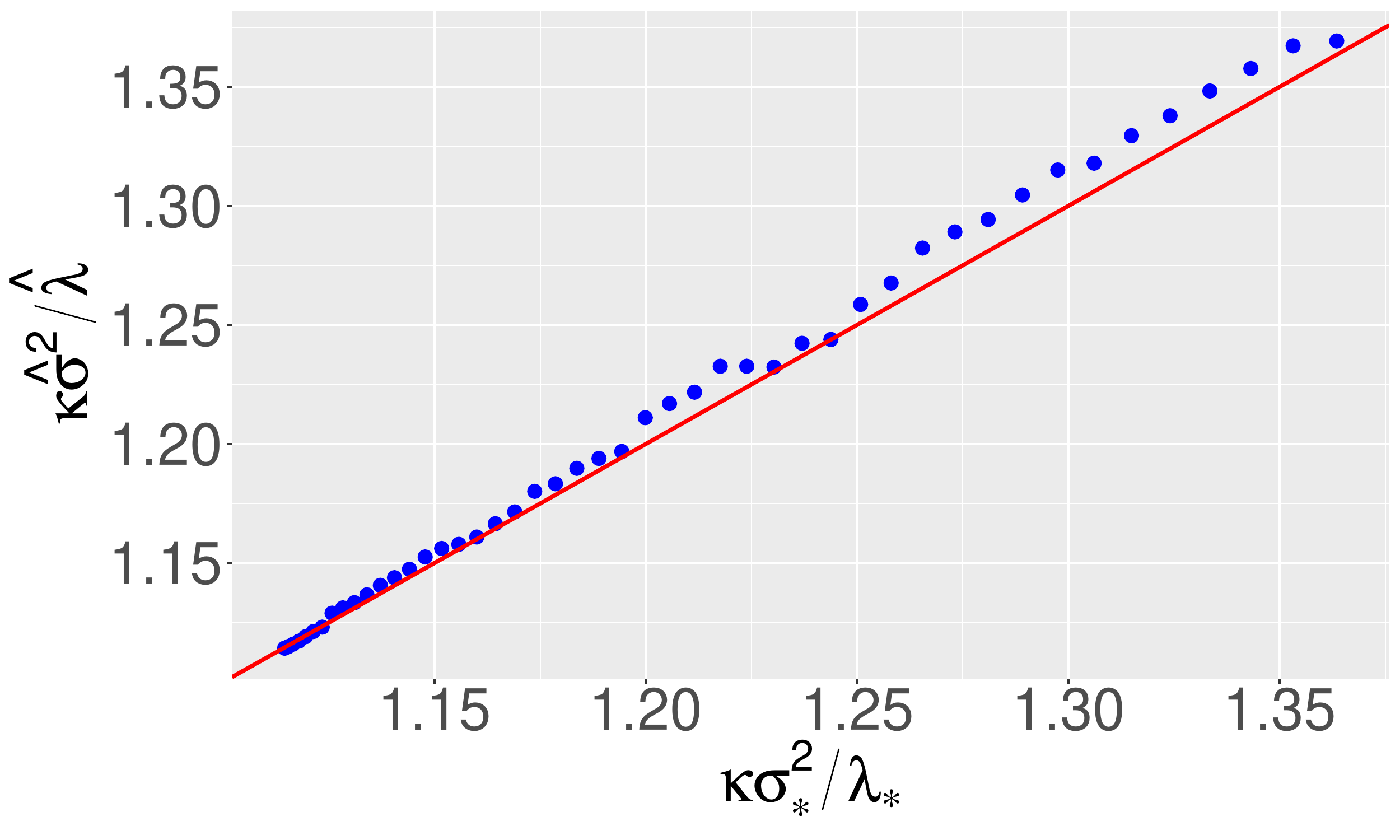}
            \caption{}    
        \end{subfigure}
        \caption{ProbeFrontier estimates of signal strength $\hat{\gamma}$, bias
          $\hat{\alpha}$, std.~dev.~$\hat{\sigma}$, and LRT factor
          $\kappa \hat{\sigma}^2/\hat \lambda$ in \eqref{eq:lrt}, 
          plotted against the theoretical values. }
       \label{fig:estimates}
    \end{figure}
    Next, we study the accuracy of $\hat{\gamma}$ across different
    choices for $\gamma$, ranging from $0.3$ to $5$. We begin by
    selecting a fine grid of $\gamma$ values and for each, we generate
    observations $(y_i,\bX_i)$ with $n=4000$, $p = 400$ (so that
    $\kappa=0.1$), and half the coefficients have a nonvanishing
    magnitude scaled in such a way that the signal strength is
    $\gamma$. Figure \ref{fig:estimates}(a) displays $\hat{\gamma}$
    versus $\gamma$ in blue, and we notice that ProbeFrontier works
    very well. We observe that the blue points fluctuate very mildly
    above the diagonal for larger values of the signal strength but
    remain extremely close to the diagonal throughout.  This confirms
    that ProbeFrontier estimates the signal strength $\gamma$ with
    reasonably high precision.  Having obtained an accurate estimate
    for $\gamma$, plugging it into \eqref{eq:main} immediately yields
    an estimate for the bias $\alphas$, standard deviation $\sigmas$
    and the rescaling factor in \eqref{eq:lrt}. We study the accuracy
    of these estimates in Figure \ref{fig:estimates}(b)-(d). We
    observe a similar behavior in all these cases, with the procedure
    yielding extremely precise estimates for smaller values, and
    reasonably accurate estimates for higher values.

\begin{figure}
\begin{center}
	\begin{tabular}{cc}
	\includegraphics[scale=0.315,keepaspectratio]{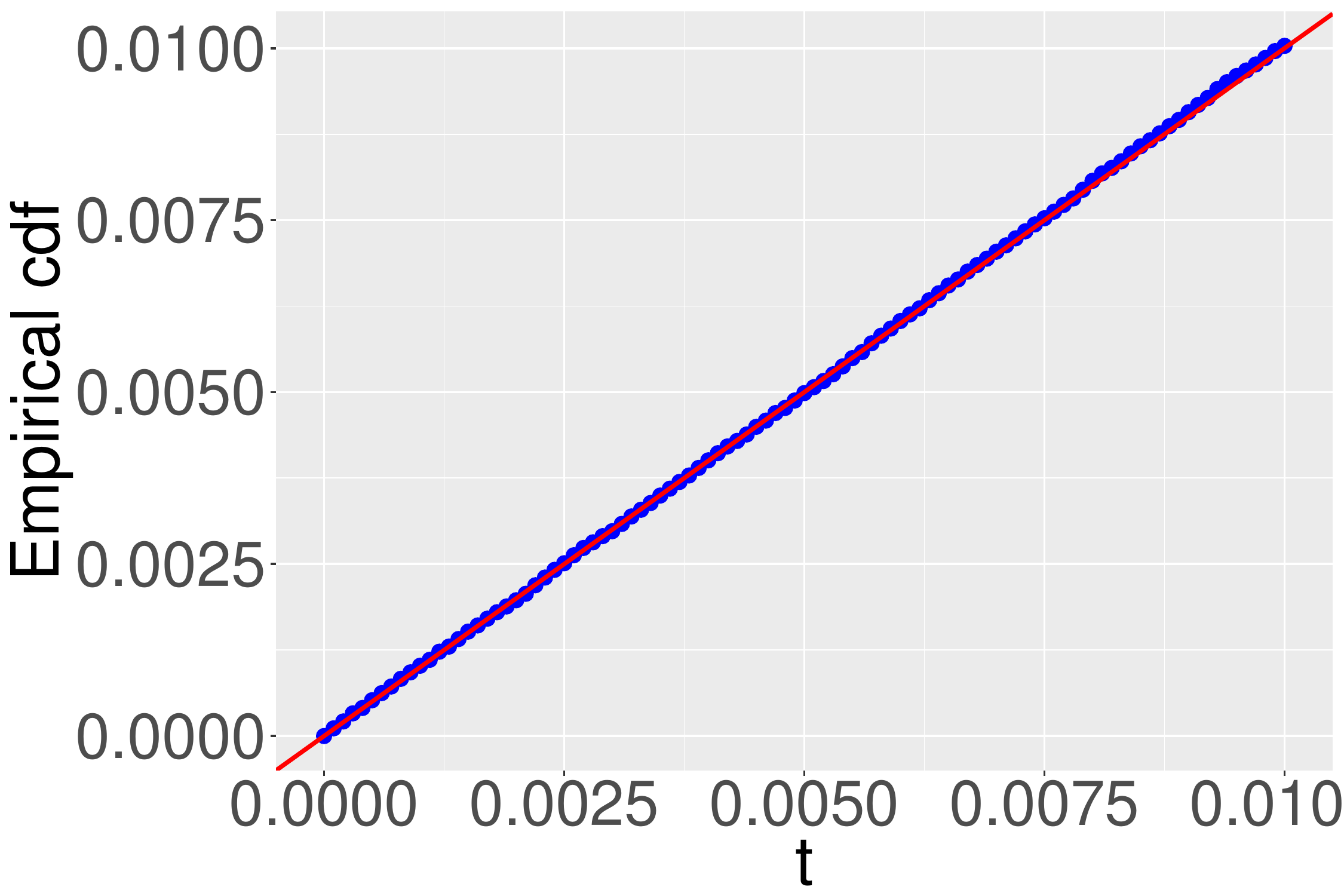} 
          &  \includegraphics[scale=0.315,keepaspectratio]{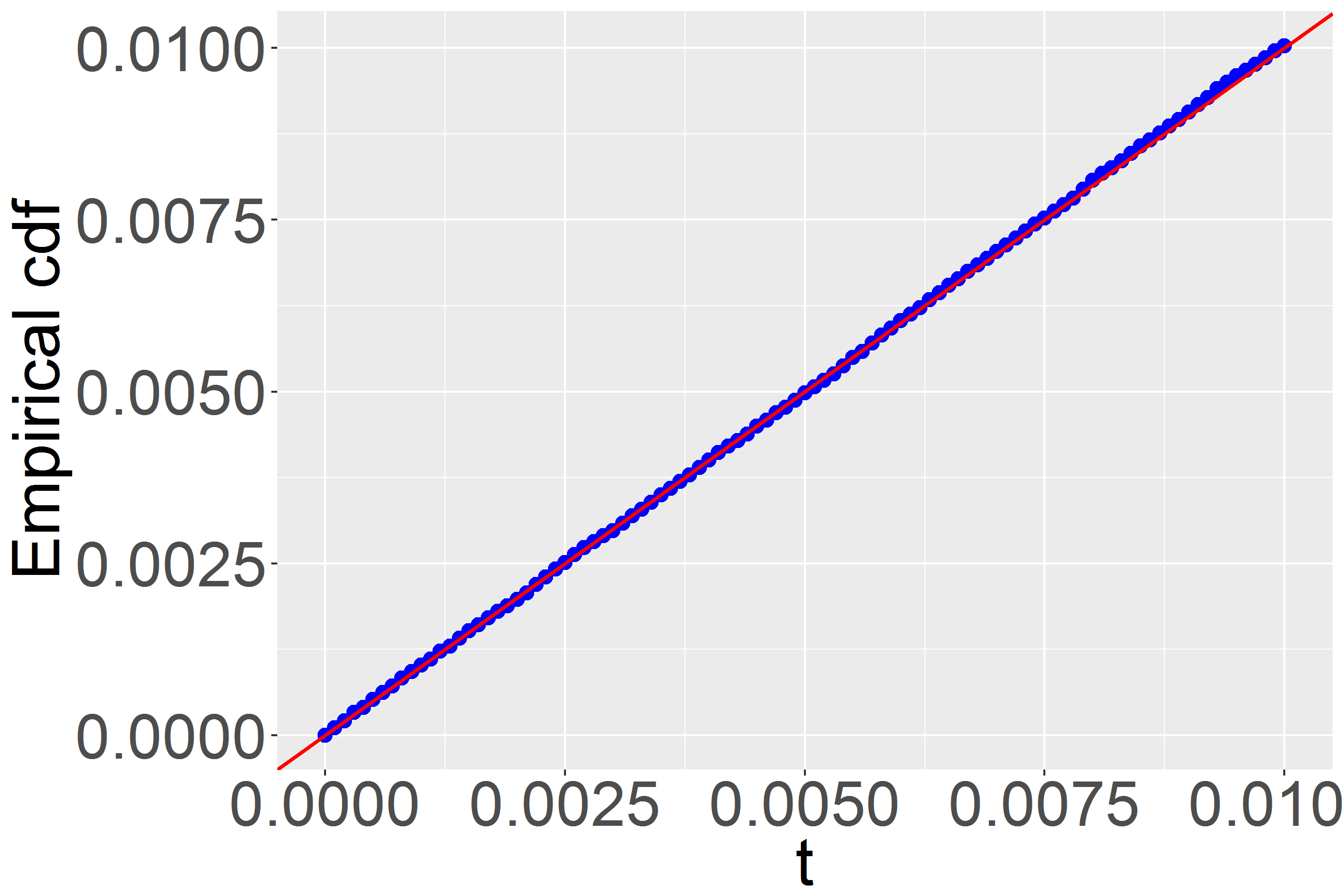}\tabularnewline
		(a) & (b)  \tabularnewline
	\end{tabular}
\end{center}
\caption{(a) Empirical distribution of the p-values based on the LLR
  approximation \eqref{eq:lrt}, obtained using the estimated factor
  $\kappa\hat{\sigma}^2/\hat{\lambda}$.  (b) Same as (a), but showing
  the tail of the empirical cdf. The calculations are based on 500,000
  replicates. 
  }
\label{fig:ecdf_est_pval}
\end{figure}

\begin{table}[h!]
\centering
\begin{tabular} {|c|c|c|c|c|}
  \hline
       Parameters & $\gamma$ & $\alpha$ & $\sigma$ & $\kappa \sigma^2/\lambda$\\
       \hline
     True & $2.2361$ & $1.1678$ & $3.3466$ & $1.166$ \\
     \hline
     Estimates & $2.2771 (0.0012)$ & $  1.1698  (0.0001)$  & $3.3751  (0.0008)$ & $1.1680 ( 0.0001)$ \\
     \hline
    \end{tabular} \\ 
    \caption{Parameter estimates in the setting of Table
      \ref{tab:alphasigma}. We report an average over 6000 replicates
      with the std.~dev.~between parentheses.}.
    \label{tab:paramest}   
  \end{table}
  Finally, we focus on the estimation accuracy for a particular
  $(\kappa,\gamma)$ pair across several replicates. In the setting of
  Figure \ref{fig:pvaladjusted}, we generate $6000$ samples and obtain
  estimates of bias ($\hat{\alpha}$), std.~dev.~($\hat{\sigma}$), and
  rescaling factor for the LRT ($\kappa \hat{\sigma}^2/\hat \lambda$).
  The average of these estimates are reported in Table
  \ref{tab:paramest}. Our estimates
  always recover the true values up to the first digit. It is
  instructive to study the precision of the procedure on the p-value
  scale. To this end, we compute p-values from \eqref{eq:lrt}, using
  the estimated multiplicative factor
  $\kappa \hat{\sigma}^2/ \hat{\lambda}$. The empirical cdf of the
  p-values both in the bulk and in the extreme tails is shown in
  Figure \ref{fig:ecdf_est_pval}. We observe the perfect agreement
  with the uniform distribution, establishing the practical
  applicability of our theory and methods.

  \subsection{De-biasing the MLE and its predictions}
  
\begin{figure}
\begin{center}
	\begin{tabular}{ccc}
	\includegraphics[scale=0.3,keepaspectratio]{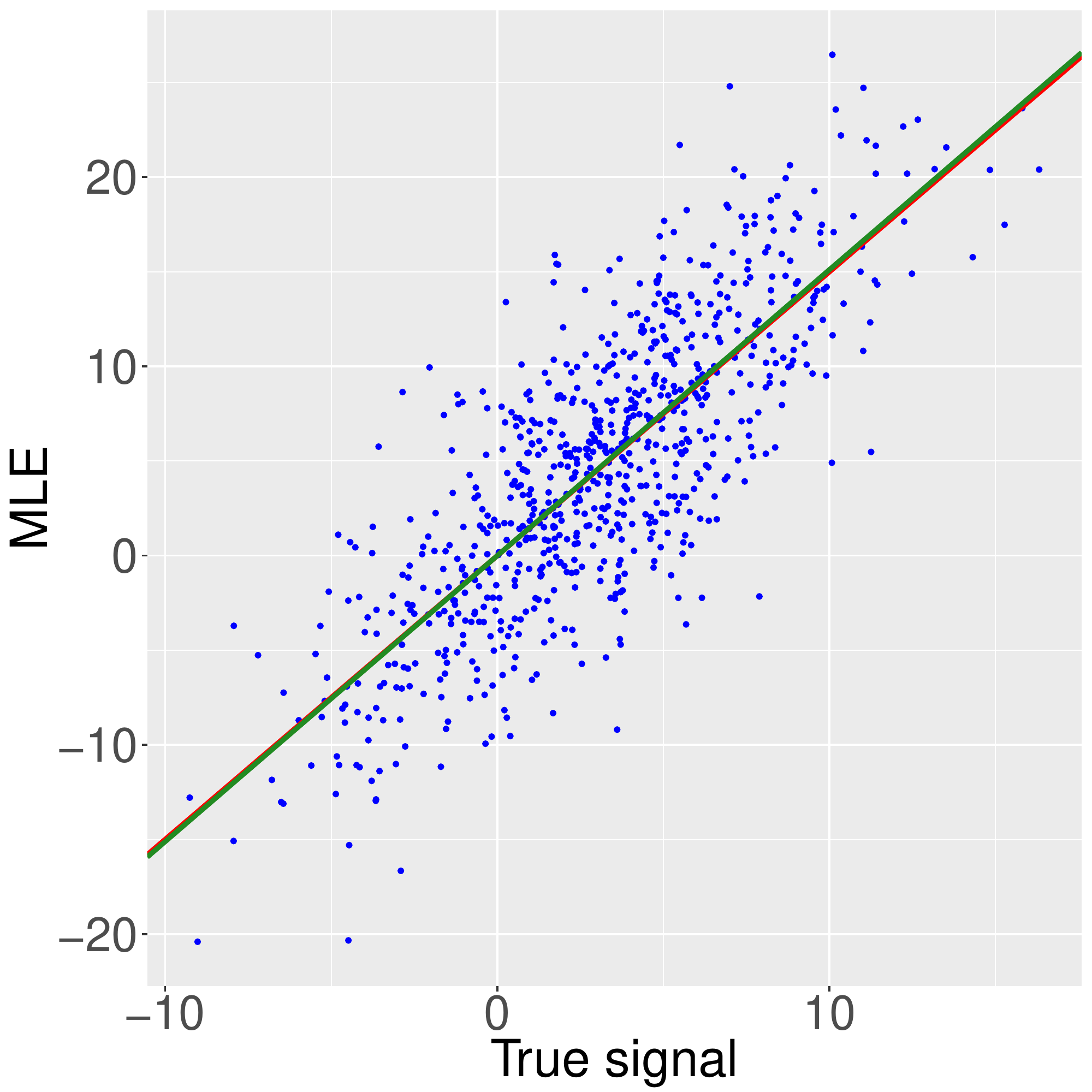} 
        & & \raisebox{0.4cm} {\includegraphics[scale=0.281,keepaspectratio]{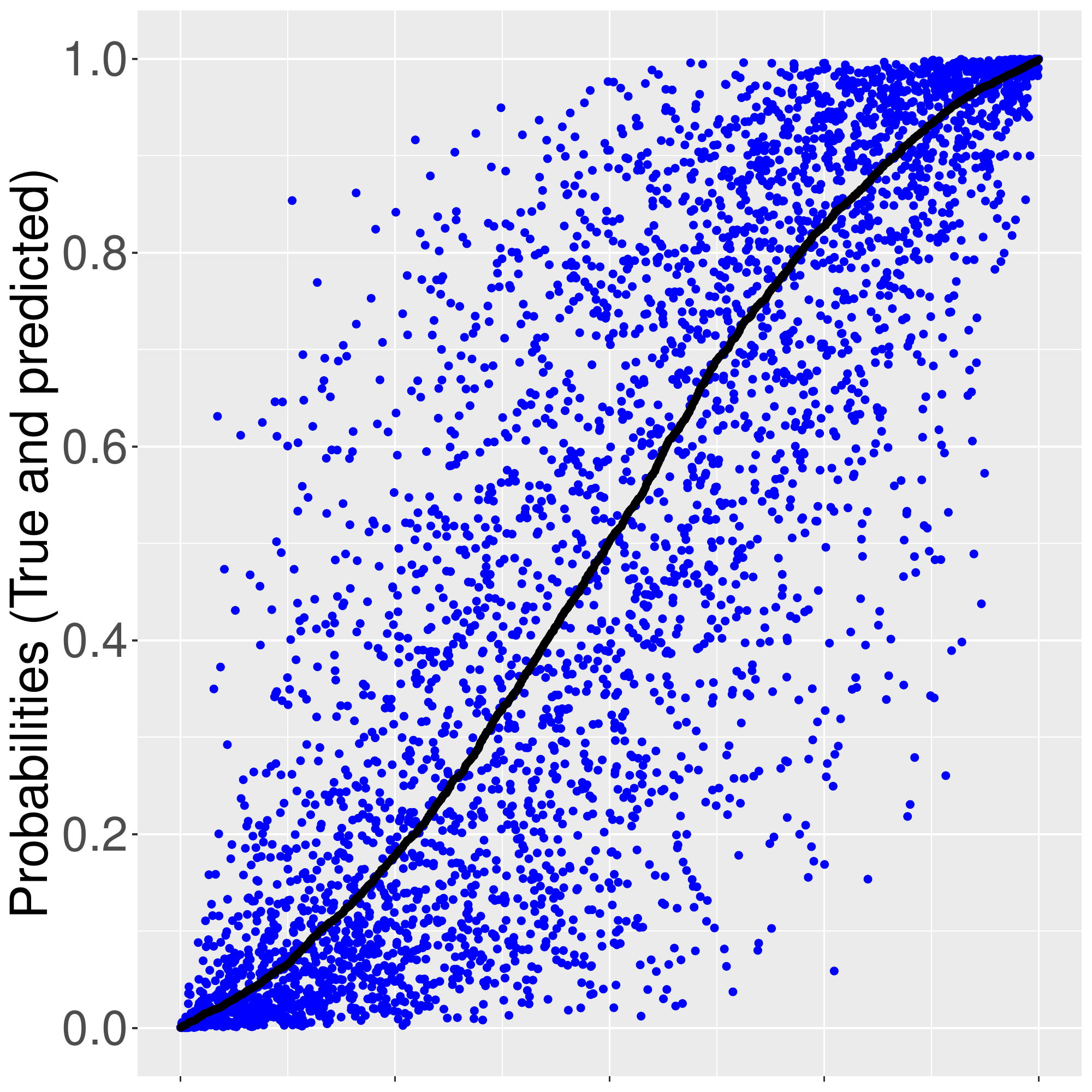}}\tabularnewline
		(a)& & (b)  \tabularnewline
	\end{tabular}
\end{center}
\caption{(a) Scatterplot of the pairs $(\beta_j, \hat{\beta}_j)$
    for the dataset from Figure \ref{fig:centering}. Here,
    $\alphas = 1.499$ (red line) and our ProbeFrontier estimate is
    $ \hat{\alpha} = 1.511$ (green line). The estimate is so close
    that the green line masks the red. (b) True conditional
    probabilities $\rho'(\bX_i' \bbeta)$ (black curve), and
    corresponding estimated probabilities
    $\rho'(\bX_i \hbbeta/\hat{\alpha})$ computed from the {\em
      de-biased MLE} (blue point). Observe that the black curve now
    passes through the center of the blue point cloud. Our predictions
    are fairly unbiased.}\label{fig:centering_est}
\end{figure}

  We have seen that maximum likelihood produces biased
    coefficient estimates and predictions.  The question is how
    precisely can our proposed theory and methods correct this. Recall
    the example from Figure \ref{fig:centering}, where the theoretical
    prediction for the bias is $\alphas = 1.499$. For this dataset,
    ProbeFrontier yields $\hat \alpha = 1.511$, shown as the green
    line in Figure
    \ref{fig:centering_est}(a). 
    Clearly, the estimate of bias is extremely precise and
    coefficient estimates $\hbeta_j/\hat \alpha$ appear nearly
    unbiased.

Further, we can also use our estimate of bias to de-bias the
  predictions since we can estimate the regression function by
  $\rho'(\bX'\hbbeta/\hat{\alpha})$. Figure \ref{fig:centering_est}(b)
  shows our predictions on the same dataset. In stark contrast to
  Figure \ref{fig:centering}(b), the predictions are now centered
  around the regression function (the method seems fairly unbiased),
  and the massive shrinkage towards the extremes has disappeared. The
  predictions constructed from the debiased MLE no longer falsely
  predict almost certain outcomes. Rather, we obtain fairly
  non-trivial chances of being classified in either of the two
  response categories---as it should be.


\section{Main Mathematical Ideas}

As we mentioned earlier, we do not provide detailed proofs in
this paper. The reader will find them in \cite{sur2018modernproofs} and the first author's Ph.~D.~thesis. 
However, in this section we give some
of the key mathematical ideas and explain some of the main steps in
the arguments, relying as much as possible on published results from
\cite{sur2017likelihood}.

\subsection{The bulk distribution of the MLE}

To analyze the MLE, we introduce an approximate message passing (AMP)
algorithm that tracks the MLE in the limit of large $n$ and $p$. Our
purpose is a little different from the work in
\cite{rangan2011generalized} which, in the context of generalized
linear models, proposed AMP algorithms for Bayesian posterior
inference, and whose properties have later been studied in
\cite{javanmard2013state} and \cite{barbier2017phase}.  To the best of
our knowledge, an AMP algorithm for tracking the MLE from a logistic
model has not yet been proposed in the literature.  Our starting point
is to write down a sequence of AMP iterates
$\{\bS^t,\hbbeta^t \}_{t \geq 0}$, with
$\bS^t \in \R^n, \hbbeta^t \in \R^p$, using the following scheme:
starting with an initial guess $\bbeta^0$, set $\bS^0 = \bX \bbeta^0$
and for $t = 1, 2, \ldots$, inductively define
\begin{equation}
 \begin{aligned}
\hbbeta^{t} & = \hbbeta^{t-1} + \kappa^{-1} \bX'\Psi_{t-1}(\by,\bS^{t-1})\\
\bS^t & = \bX \hbbeta^t - \Psi_{t-1} (\by,\bS^{t-1}) 
\end{aligned}
\label{eq:AMP}
\end{equation}
where the function $\Psi_t $ is applied element-wise and is equal to 
\begin{equation}\label{eq:Psi_t}
  \Psi_t(y,s) = \lambda_t r_t, \qquad r_t = y - \rho'(\prox_{\lambda_t \rho} (\lambda_t y
  +s)). 
\end{equation}
 Observe that the evolution \eqref{eq:AMP} depends on
a sequence of parameters $\{\lambda_t\}$ whose dynamics we describe
next.

This description requires introducing an augmented sequence
$\{\alpha_t, \sigma_t, \lambda_t\}_{t \ge 0}$. With these two extra
parameters $(\alpha_t, \sigma_t)$, we let $(Q_1^t, Q_2^t)$ be a
bivariate normal variable with mean $\bm{0}$ and covariance matrix
$\bSigma( \alpha_t, \sigma_t)$ defined exactly as in
\eqref{eq:covfunc}. Then starting from an initial pair
$\alpha_0, \sigma_0$, for $t = 0, 1, \ldots$, we inductively define
$\lambda_t$ as the solution to
\begin{equation}
  \label{eq:varmap}
  \E\left[\frac{2\rho'(Q_1^t)}{1+\lambda \rho''(\prox_{\lambda \rho}(Q_2^t))} \right]  = 1- \kappa
\end{equation}
and the extra parameters $\alpha_{t+1}, \sigma_{t+1}$ as
  \begin{equation} 
\begin{aligned} 
  \alpha_{t+1} & = \alpha_t + \frac{1}{\kappa \gamma^2}\E \left[2\rho'(Q_1^t)Q_1^t\lambda_t \rho'(\prox_{\lambda_t \rho}(Q_2^t)) \right]  \\
  \sigma_{t+1}^2 & = \frac{1}{\kappa^2}\E \left[ 2\rho'(Q_1^t) \left(\lambda_t \rho'(\prox_{\lambda_t \rho} (Q_2^t)) \right)^2\right].
  \end{aligned}
\label{eq:varmaptwo}
\end{equation}
To repeat, we run the AMP iterations \eqref{eq:AMP} using the scalar
variables $\{\lambda_t\}$ calculated via the {\em variance map}
updates \eqref{eq:varmap}--\eqref{eq:varmaptwo}. 

In the regime where the MLE exists (see Figure \ref{fig:kappagamma}),
the variance map updates \eqref{eq:varmap}--\eqref{eq:varmaptwo}
converge (as $t \rightarrow \infty$) to a unique fixed point $(\alphas,\sigmas,\lambdas)$. Note
that by definition, $(\alphas, \sigmas, \lambdas)$ is the solution to
our system \eqref{eq:main} in three unknowns. From now on, we use 
$\alpha_0 = \alphas$, $\sigma_0 = \sigmas$ so that the sequence
$\{\alpha_t, \sigma_t, \lambda_t\}$ is stationary; i.~e.~for all $t \ge 0$,
\[
  \alpha_t = \alphas, \quad \sigma_t = \sigmas, \quad \lambda_t = \lambdas.
\]
With this stationary sequence of parameters, imagine now initializing
the AMP iterations with a vector $\hbbeta^0$ obeying
  $$\lim_{n,p \rightarrow \infty} \, \frac{1}{p} \|\hbbeta^0 - \alphas \bbeta \|^2 =
  \sigmas^2.$$ It is not hard to see that if the proposed AMP
  algorithm converges to a fixed point
  $\{\bS_{\star},\hbbeta_{\star} \}$, then it is such that
  $\nabla \ell(\hbbeta_{\star})=0$ (see Appendix
  \ref{sec:optimality}); that is, $\hbbeta_{\star}$ obeys the MLE
  optimality conditions.  This provides some intuition as to why the
  above algorithm would turn out to be useful in this context.
  
  
    The crucial point is that we can
  study the properties of the MLE by studying the properties of the
  AMP iterates with the proviso that they converge. It turns out that
  the study of the sequence $\{\bS^t, \hbbeta^t\}$ is {amenable to a
    rigorous analysis because several transformations reduce the above
    algorithm to a generalized AMP algorithm
    \cite{javanmard2013state}, which in turn yields a characterization
    of the limiting variance of the AMP iterates:} for any function
  $\psi$ as in Theorem \ref{thm:mle}, we have {as $n \rightarrow \infty$}, 
  \begin{equation}\label{eq:stateev}
  \frac{1}{p} \sum_{j=1}^{p}\psi(\hat{\beta}_{j}^{t} - \alphas \beta_{j},\beta_j)   \,\, {\stackrel{\text{a.s.}}{\longrightarrow}} \, \, \E \left[ \psi(\sigmas Z,\beta) \right],
  \end{equation}
  where $\beta$ is drawn {from the distribution $\Pi$} (see Theorem
  \ref{thm:mle}) independently of $Z\sim\dnorm(0,1)$, and $\sigmas$ is
  as above. {To summarize, the asymptotic behavior of the AMP iterates
    $\hbbeta^t$ can be characterized through a standard Gaussian
    variable, the distribution $\Pi$ and the scalar quantity $\sigmas$
    determined by the iteration
    \eqref{eq:varmap}--\eqref{eq:varmaptwo}. The description of our
    AMP algorithm and large sample properties of the iterates are
    understood only when we understand the behavior of the scalar
    sequences $\{\alpha_t, \sigma_t, \lambda_t \}_{t \geq 0}$, which
    are known as the state evolution sequence in the literature; this
    formalism was introduced in
    \cite{bayati2011dynamics,donoho2009message,donoho2010message,donoho2011noise}.}
  From here on, an analysis similar to that in
  \cite{sur2017likelihood} establishes that in the limit of large
  iteration counts, the AMP iterates converge to the MLE, that is,
\[\lim_{t \rightarrow \infty} \lim_{n \rightarrow \infty}\frac{1}{p} \sum_{j=1}^{p}\psi(\hbeta_{j}^{t } - \alphas \beta_{j}, \beta_j) = \lim_{n \rightarrow \infty}\frac{1}{p} \sum_{j=1}^{p}\psi(\hbeta_{j} - \alphas \beta_{j}, \beta_j),\]
which is the content of Theorem \ref{thm:mle}.  

\subsection{The distribution of a null coordinate}\label{subsec:null}

We sketch the proof of Theorem \ref{thm:mainthm} in the case
  where the empirical limiting distribution $\Pi$ has a point mass at
  zero. The analysis in the general case, where the number of
  vanishing coefficients is arbitrary, and in particular, $o(n)$, is
  very different and may be found in Appendix \ref{sec:refined}.


Now consider Theorem \ref{thm:mle} with
  $\psi(t,u) = t^2 1(u = 0)$. Strictly speaking, this is a
  discontinuous function which is not pseudo-Lipschitz.  However, we 
  can work with a smooth approximation $\psi_a$, instead, obtained
  using standard techniques for smoothing an indicator function, such
  that the error $\|\psi -\psi_a \|_{2}$ is arbitrarily
  small. For simplicity, we skip the technical details underlying this
  approximation, and motivate the subsequent arguments using $\psi$
  directly. Theorem \ref{thm:mle} then yields
\begin{equation}\label{eq:norm}
\frac{1}{p} \sum_{j \in [p]: \beta_j =0} \hat{\beta}_j^2  \,\, {\stackrel{\text{a.s.}}{\longrightarrow}}\,\, \sigmas^2 \prob_{\Pi}\left[\beta =0  \right] \ \ \implies \frac{1}{|j\in [p]:\beta_j = 0|} \sum_{j \in [p]: \beta_j =0} \hat{\beta}_j^2 \,\,  {\stackrel{\text{a.s.}}{\longrightarrow}} \,\, \sigmas^2. 
\end{equation}
Without loss of generality, assume that the first $k$ coordinates of
$\bbeta$ vanish, and that $\bbeta$ is partitioned as
$\bbeta = \left(\bzero_{[k]},\bbeta_{-[k]} \right)$ and similarly for
$\hbbeta$. From the rotational distributional invariance of the
$\bX_i$'s, it can be shown that for any fixed orthogonal matrix
$\bU \in \R^{k \times k}$,
$\hbbeta \eqd \left(\bU \hbbeta_{[k]},
  \hbbeta_{-[k]}\right)$. Consequently,
$\hbbeta_{[k]} / \|\hbbeta_{[k]} \|$ is uniformly distributed on the
unit sphere $\mathbb{S}^{k-1}$ and is independent of
$\|\hbbeta_{[k]} \|$. Thus, any null coordinate $\hat{\beta}_j$ has
the same distribution as $\|\hbbeta_{[k]}\| \, Z_j /\|\bZ \|$, where
$\bZ \sim \dnorm(0, \bm{I}_k)$, independent of $\hbbeta_{[k]}$. From
\eqref{eq:norm} and the weak law of large numbers, we have
$\|\hbbeta_{[k]} \|/\|\bZ \| \convP \sigmas$, leading to
$\hat{\beta}_j \convd \dnorm(0, \sigmas^2)$.

\subsection{The distribution of the LRT}\label{subsec:lrtdist}

Once the distribution of $\hbeta_j$ for a null $j$ is known, the
distribution of the LRT is a stone throw away, at least conceptually;
that is to say, if we are willing to ignore some technical
difficulties and leverage existing work. Indeed, following a reduction
similar to that in \cite{sur2017likelihood}, one can establish that
\begin{equation}\label{eq:lrtapprox}
2\Lambda_{j} = \frac{\kappa}{{\lambda_{[-j]}}} \hbeta_j^2 + o_P(1), 
\end{equation}
where
$\lambda_{[-j]} := \tr \left[\nabla^2(\ell_{[-j]}(\hbbeta_{[-j]}))^{-1}\right]/n$
in which $\ell_{[-j]}$ is the negative log-likelihood with the $j$-th
variable removed and $\hbbeta_{[-j]}$ the corresponding MLE.
Put $\lambda = \tr [\nabla^2(\ell(\hbbeta))^{-1}]/n$. Then following an
approach similar to that in \cite[Appendix I]{sur2017likelihood}, it
can be established that
$\lambda_{[-j]} = \lambda + o_P(1) \convP \lambdas$.  This gives that
$2\Lambda_{j}$ is a multiple of a $\chi_1^2$ variable with
multiplicative factor given by $\kappa \sigmas^2/\lambdas$.

This rough analysis shows that the distribution of the LLR in high
dimensions deviates from a $\chi^2_1$ due to the coupled effects of
two high-dimensional phenomena. The first is the inflated variance of
the MLE, which is larger than classically predicted. The second comes
from the term $\lambdas$, which is approximately equal to
$\tr\left(\bm{H}^{-1}(\hbbeta)\right)/n$, where
$\bm{H}(\hbbeta) = \nabla^2 \ell(\hbbeta)$ is the Hessian of the
negative log-likelihood function. In the classical setting, this
Hessian converges to a population limit. This is not the case in
higher dimensions and the greater spread in the eigenvalues also
contributes to the magnitude of the LRT.

\section{Broader Implications and Future Directions}

This paper shows that in high-dimensions, classical ML theory is
unacceptable. Among other things, classical theory predicts that the
MLE is approximately unbiased when in reality it seriously
overestimates effect magnitudes. Since the purpose of logistic
modeling is to estimate the risk of a specific disease given a
patient's observed characteristics, say, the bias of the MLE is
extremely problematic. As we have seen, an immediate consequence of
the strong bias is that the MLE either dramatically overestimates, or
underestimates, the chance of being sick. The issue becomes
increasingly severe as either the dimensionality or the signal
strength, or both, increase. This, along with the fact that p-values
computed from classical approximations are misleading, clearly make
the case that routinely used statistical tools fail to provide
meaningful inferences from both an estimation and a testing
perspective.

We have developed a new theory which gives the asymptotic distribution
of the MLE and the LRT in a model with independent covariates. As seen
in Section \ref{sec:nonGaussian}, our results likely hold for a broader
range of feature distributions (i.e.~other than Gaussian) and it would
be important to establish this rigorously.  Further, we have also
shown how to adjust inference by plugging in an estimate of signal
strength in our theoretical predictions. Although our method works
extremely well, it would be of interest to study others as well.

Finally, we conclude this paper with two promising directions for
future work: (1) It would be of great interest to develop
corresponding results in the case where the predictors are correlated.
(2) It would be of interest to extend the results from this paper to
other generalized linear models.

\subsection*{Acknowledgements}

P.~S.~was partially supported by the Ric Weiland Graduate Fellowship in
the School of Humanities and Sciences, Stanford University. E.~C.~was
partially supported by the Office of Naval Research under grant
N00014-16-1-2712, by the National Science Foundation via DMS 1712800,
by the Math + X Award from the Simons Foundation and by a generous
gift from TwoSigma. P.~S.~would like to thank Andrea Montanari and
Subhabrata Sen for helpful discussions on AMP. We thank Ma\l{}gorzata
Bogdan, Nikolaos Ignatiadis, Asaf Weinstein, Lucas Janson and Chiara Sabatti for useful comments
about an early version of the paper. 

\bibliographystyle{plain}
\bibliography{bibfileLR}

\normalsize 
\newpage 
\appendix 
\section{Fisher information}
\label{sec:Fisher}

We work with the model from Section \ref{sec:main} and introduce the
Fisher information matrix defined as
\[
I(\bbeta) = \E \left[ \sum_i \rho''(\bX_i'\bbeta) \bX_i \bX_i'\right] = n  \E \left[ \rho''(\bX_i'\bbeta) \bX_i \bX_i'\right].  
\]
With $\bX_i \sim \dnorm(\bzero, n^{-1}{\bm I})$, it is not hard to see
that the $(k,j)$th entry of the matrix
$n \rho''(\bX_i'\bbeta) \bX_i \bX_i'$ is distributed as
\[
 \rho''(\gamma X_1) X_k X_j, \quad X_1, \ldots, X_p \iid \dnorm(0,1). 
\]
From here on, a reasonably straightforward calculation gives
\[
I(\bbeta) = \nu ({\bm I} + \delta \bu \bu'),\quad \bu = \bbeta/\|\bbeta\|, 
\]
where 
\[
  \nu = \E [\rho''(\gamma X_1)], \qquad \delta = \frac{\E
    [\rho''(\gamma X_1) X_1^2] - \E [\rho''(\gamma X_1)]}{\E [\rho''(\gamma X_1)]}.
\]
This implies that 
\[
I^{-1}(\bbeta) = \nu^{-1} \left({\bm I} - \frac{\delta}{1+\delta} \bu\bu'\right),
\]
which means that the classically predicted variance of $\hat{\beta}_j$
is equal to
\[
  \nu^{-1} \left(1 - \frac{\delta}{1+\delta}
    \frac{\beta_j^2}{\|\bbeta\|^2}\right).
\]
When $\beta_j = 0$, the predicted standard deviation is
$\nu^{-1/2} = 2.66$ for $\gamma^2 = 5$.

Statistical software packages base their inferences on the approximate
Fisher information defined as
$\sum_i \rho''(\bX_i'\hbbeta) \bX_i \bX_i'$ (or small corrections
thereof). This treats the covariates as fixed and substitutes the
value of the unknown regression coefficient sequence $\bbeta$ with
that of the MLE $\hbbeta$ (plugin estimate).

\section{Properties of fixed points of the AMP algorithm}
\label{sec:optimality}

In this section, we elaborate on the connection between the fixed points of \eqref{eq:AMP} and the MLE $\hbbeta$. From \eqref{eq:AMP}, we immediately see that if $(\hbbeta_{\star},\bS_{\star})$ is a fixed point, the pair satisfies
\begin{align*}
 \bX'\{ \by - \rho'(\prox_{\lambdas \rho }(\lambdas \by+\bS_{\star})) \} & =\bzero  \\
 (\lambdas \by +\bS_{\star} ) - \lambdas \rho'(\prox_{\lambdas \rho} (\lambdas \by+\bS_{\star})) & = \bX \hbbeta_{\star}.
 \end{align*}
 Since by definition of the proximal mapping operator,
 $ z-\lambda \rho'(\prox_{\lambda \rho }(z) ) = \prox_{\lambda \rho
 }(z)$, we have that
 $\bX \hbbeta_{\star} = \prox_{\lambdas \rho} (\lambdas
 \by+\bS_{\star}) $ which implies
 \[ \bX'\{ \by - \rho'(\bX \hbbeta_{\star}) \} = \bzero. \] Hence, the
 fixed point $\hbbeta_{\star}$ obeys
 $\nabla \ell(\hbbeta_{\star}) = \bzero$, the optimality condition for the MLE.
 
 \section{Refined analysis of the distribution of a null coordinate}
 \label{sec:refined}
The AMP analysis is useful to analyze the bulk behavior of the MLE;
i.e.~the expected behavior when averaging over all coordinates. It
also helps in characterizing the distribution of a null coordinate
when the limiting empirical cdf does not have a point mass at zero, as
we have seen in Section \ref{subsec:null}.  However, the study of the
behavior of a single coordinate when there is an arbitrary number of
nulls requires a more refined analysis. To this end, the proof uses a
leave-one-out approach, as in \cite{el2015impact,el2013robust,
  sur2017likelihood}. The complete rigorous technical details are very
involved and this is a reason why we only present approximate or
non-rigorous heuristic calculations.

Fix $j$ such that $\beta_j = 0$. 
Since the corresponding
  predictor plays no role in the distribution of the response, we
  expect that including this predictor or not in the regression will
  not make much difference in the fitted values, that is,
\begin{equation}
  \label{eq:approx}
\bX_i'\hbbeta \approx \bX_{i,-j}' \hbbeta_{[-j]};
\end{equation}
here, $\bX_{i,-j}$ is $i$-th row of the reduced matrix of predictors
with the $j$-th column removed and $\hbbeta_{[-j]}$ is the MLE for the
reduced model.
On the one hand, the approximation \eqref{eq:approx} suggests Taylor
expanding $\rho'(\bX_i'\hbbeta)$ around the point
$\bX_{i,-j}'\hbbeta_{[-j]}$:
\[\rho'(\bX_i'\hbbeta) \approx  \rho'(\bX_{i,-j}'\hbbeta_{[-j]}) + \rho''(\bX_{i,-j}'\hbbeta_{[-j]}) \left[X_{ij} \hat{\beta}_{j} + \bX_{i,-j}'\left(\hbbeta_{-j} - \hbbeta_{[-j]} \right) \right], \]
where $\hbbeta_{-j}$ is the full-model MLE vector, however, without
the $j$-th coordinate.  On the other hand, we can subtract the
two score equations $\nabla \ell (\hbbeta)=0$ and
$ \nabla {\ell_{[-j]}} (\hbbeta_{[-j]})=0$ ($\ell_{[-j]}$ is the
negative log-likelihood for the reduced model), which upon separating
the components corresponding to the $j$-th coordinate from the others,
yields
\begin{align*}
\sum_{i=1}^{n} X_{ij} \left(y_i - \rho'(\bX_{i}' \hbbeta) \right) & = 0\\
 \sum_{i=1}^{n} \bX_{i,-j} \{\rho'(\bX_{i}' \hbbeta) - \rho'(\bX_{i,-j}'\hbbeta_{[-j]}) \} & = \bzero. 
\end{align*}  
Plugging in the approximation for $\rho'(\bX_{i}' \hbbeta)$ yields two
equations in the two unknowns $\hat{\beta}_{j}$ and
$(\hbbeta_{-j} - \hbbeta_{[-j]})$. After some algebra, solving for
$\hat{\beta}_{j}$ yields 
\begin{equation*}
\hat{\beta}_j = \frac{\sum_{i=1}^{n}X_{ij}\left(y_i - \rho'(\bX_{i,-j}'\hbbeta_{[-j]}) \right)}{\bX_{\bullet j}'\bD^{1/2}(\hbbeta_{[-j]})\bH \bD^{1/2}(\hbbeta_{[-j]})\bX_{\bullet j}} + o_P(1),
\end{equation*}
where
$\bH = \bm{I} - \bD^{1/2}(\hbbeta_{[-j]})\bX_{\bullet -j}(\nabla^2
\ell_{-j}(\hbbeta_{[-j]}) )^{-1}\bX_{\bullet
  -j}'\bD^{1/2}(\hbbeta_{[-j]})$ and $\bD(\hbbeta_{[-j]})$ is an
$n \times n$ diagonal matrix with $i-$th entry given by
$\rho''(\bX_{i,-j}'\hbbeta_{[-j]})$.  Above $\bX_{\bullet j}$ is the
$j$-th column of $\bX$ and $\bX_{\bullet -j}$ all the others.  It was established in \cite{sur2017likelihood} that the denominator above is equal to $\kappa /\lambda_{[-j]} + o_P(1)$, {where, we have see in Section \ref{subsec:lrtdist} that
$$
\lambda_{[-j]} := \frac{1}{n} \tr [\nabla^2(\ell_{[-j]}(\hbbeta_{[-j]}))^{-1}].$$ }
Note that since $\beta_j = 0$, $\by$ and $\bX_{\bullet -j}, \hbbeta_{[-j]}$
are independent of $\bX_{\bullet j}$. This gives the approximation 
\begin{equation}\label{eq:betahatj}
  \hat{\beta}_j = \frac{\lambda_{[-j]} s_j}{\kappa}  \, Z + o_P(1), \quad s_j^2 = \frac{1}{n}\sum_{i=1}^n\left(y_i
    - \rho'(\bX_{i,-j}' \hbbeta_{[-j]})\right)^2,
\end{equation}
where $Z$ is an independent standard normal. In Section
  \ref{subsec:lrtdist}, we saw that $\lambda_{[-j]} \convP \lambdas$.
  It remains to understand the behavior of $s_j$. Looking at $s_j$,
the complicated dependence structure between $\hbbeta$ and $(\by,\bX)$
makes this a potentially hard task. This is why we shall use a
leave-one-out argument and seek to express the fitted values
$\bX_{i,-j}'\hbbeta_{[-j]}$ in terms of
$\bX_{i,-j}'\hbbeta_{[-i],[-j]}$, where $\hbbeta_{[-i],[-j]}$ is the
MLE when both the $j$-th predictor and the $i$-th observation are
dropped. The independence between $\bX_{i,-j}$ and
$\hbbeta_{[-i],[-j]}$ will simplify matters.  Denote by
$\nabla {\ell}_{[-i],[-j]}(\tbbeta_{[-i],[-j]})=0$ the reduced score
equation and subtract it from the score equation for $\hbbeta$ to
obtain
\[\bX_{i,-j}\left(y_i - \rho'(\bX_{i,-j}' \hbbeta_{[-j]} ) \right) + \sum_{k \neq i} \bX_{k,-j}\left( \rho'(\bX_{k,-j}' \hbbeta_{[-i],[-j]}) - \rho'(\bX_{k,-j}' \hbbeta_{[-j]})\right) = \bzero.  \]
We argue that since the number of observations is large and the
observations are i.i.d., dropping one observation is not expected to
create much of a difference in the fitted values, hence
$\bX_{k,-j}' \hbbeta_{[-i],[-j]} \approx \bX_{k,-j}' \hbbeta_{[-j]}$. A
Taylor expansion of $\rho'(\bX_{k,-j}' \hbbeta_{[-j]})$ around
the point $\bX_{k,-j}' \hbbeta_{[-i],[-j]}$ yields 
\[\bX_{i,-j}'\left(\hbbeta_{[-j]} -\hbbeta_{[-i],[-j]}\right) \approx \bX_{i,-j}'\left[\nabla^2 \ell_{[-i],[-j]}(\hbbeta_{[-i],[-j]})\right]^{-1} \bX_{i,-j} \left(y_i - \rho'(\bX_{i,-j}'\hbbeta_{[-j]}) \right). \]
Since $\bX_{i,-j}$ and $\hbbeta_{[-i],[-j]}$ are independent, by
Hanson-Wright inequality \cite[Theorem $1.1$]{rudelson2013hanson}, the
quadratic form above is approximately equal to 
$\tr\left[\nabla^2 {\ell}_{[-i],[-j]}(\hbbeta_{[-i],[-j]})^{-1}
\right]$. Recall that
$\lambda_{[-j]} = \tr[\nabla^{2}\ell_{[-j]}(\hbbeta_{[-j]})^{-1}]$ and
again, for a large number of i.i.d.~observations, we expect these two
quantities to be close. Hence, the fitted values can be approximated
as
\[\bX_{i,-j}'\hbbeta_{[-j]} \approx \bX_{i,-j}' \hbbeta_{[-i],[-j]} + \lambda_{[-j]}  \left(y_i - \rho'(\bX_{i,-j}'\hbbeta_{[-j]}) \right).\]
Recalling the definition of the proximal mapping operator, $\prox_{\lambda \rho}(z) + \lambda\rho'(\prox_{\lambda \rho}(z)) = z$, note that the above relation gives a useful approximation for the fitted values, 
\[\bX_{i,-j}'\hbbeta_{[-j]} \approx \prox_{\lambda_{[-j]}\rho} \left(\lambda_{[-j]} y_i + \bX_{i,-j}'\hbbeta_{[-i],[-j]}  \right).\]
Further, by the triangle inequality we can show that
$$\prox_{\lambda_{[-j]} \rho} \left(\lambda_{[-j]} y_i +
  \bX_{i,-j}'\hbbeta_{[-i],[-j]} \right) \approx \prox_{\lambdas \rho}
\left(\lambdas y_i + \bX_{i,-j}'\hbbeta_{[-i],[-j]} \right).$$ It can be
shown that the residuals
$\{y_i - \rho'(\prox_{\lambdas \rho} (\lambdas y_i +
\bX_{i,-j}'\hbbeta_{[-i],[-j]} ) ) \}_{i=1,\hdots, n}$ are
asymptotically independent among themselves,  
which implies that averaging over the
squared residuals  as in \eqref{eq:betahatj} should converge in
probability to the corresponding expectation, leading to
\[\hat{\beta}_j \convd \dnorm(0,\sigma^2), \quad \sigma^2 :=\frac{\lambdas^2}{\kappa^2} \lim_{n \rightarrow \infty} \E \left[y_i - \rho'\left(\prox_{\lambdas \rho} \left(\lambdas y_i + \bX_{i,-j}'\hbbeta_{[-i],[-j]} \right)\right)\right]^2. \]

To complete the analysis, it remains to characterize the asymptotic
joint distribution of
$\bX_{i,-j}' \hbbeta_{[-i],[-j]}$ and $\bX_{i}'\bbeta$ or, equivalently,  $\bX_{i,-j}'\bbeta_{-j}$ ($\bbeta_{-j}$ is the true signal with the j-th coordinate
removed) since $\beta_j =0$. 
Instead, we find the joint distribution of
$(\bX_{i,-j}'\bbeta_{-j},\bX_{i,-j}'(\hbbeta_{[-i],[-j]} - \alphas
\bbeta_{-j}) )$ and denote this pair as
$(Q_{1}^{\star}, Q_{2}^{\star})$. The asymptotic variance of
$Q_{1}^{\star}$ is given by $\gamma^2$, that of $Q_{2}^{\star}$ 
by $\kappa \sigmas^2$, while the asymptotic covariance is equal to 
\begin{align}
\lim_{n \rightarrow \infty} \frac{\langle \hbbeta_{[-i],[-j]}-\alphas\bbeta_{-j},\bbeta_{-j} \rangle}{n} & = \kappa \lim_{t \rightarrow \infty} \lim_{n \rightarrow \infty} \frac{\langle \hbbeta^t - \alphas \bbeta,\bbeta \rangle}{p} = 0,
\end{align}
by an application of \eqref{eq:stateev}. Hence, writing
$y_i = 1(U_i \le \rho'(\bX_{i,-j}'\bbeta_{-j}))$, where the $U_i$'s
are i.i.d.~$U(0,1)$ independent from anything else, we have 
\[
  \lim_{n \rightarrow \infty} \text{Var}(\hat{\beta}_j) =
  \frac{1}{\kappa^2}\lambdas^2 \E \left[ 1(U_i \le \rho'(Q_1^{\star})) -
    \rho'(\prox_{\lambdas \rho}(\alphas Q_1^{\star}+Q_2^{\star}+\lambdas 1(U_i \le
    \rho'(Q_1^\star))) \right]^2.
\]
Using this later fact, the above expression can be simplified to
$$
\frac{1}{\kappa^2} \E\left[2 \rho'(-Q_1^\star) \left(\lambdas
    \rho'(\prox_{\lambdas \rho} ( \alphas
    Q_1^{\star}+Q_2^{\star})\right)^2 \right].$$ Note that the joint
distribution of $(-Q_1^\star,\alphas Q_1^{\star}+Q_2^{\star})$ is precisely
the same as $\bSigma(\alphas, \sigmas)$ as specified by
\eqref{eq:covfunc}. Hence, recalling \eqref{eq:main}, we obtain the
asymptotic variance of $\hat{\beta}_j$ to be $\sigmas^2$.


\end{document}